\documentclass[a4paper,english,oneside,11pt]{memoir}
\usepackage{babel}
\usepackage[T1]{fontenc}
\usepackage[utf8]{inputenc}
\usepackage{fouriernc}
\usepackage{url}
\usepackage{hyperref}
\usepackage{amssymb}
\usepackage{amsmath}
\usepackage{mathtools}
\usepackage{bm}
\usepackage{xcolor}
\usepackage[ntheorem,footnoteinside=false]{mdframed}
\usepackage[amsmath,hyperref]{ntheorem}

\newcommand{\C}{\mathbb{C}}
\newcommand{\Z}{\mathbb{Z}}
\newcommand{\N}{\mathbb{N}}
\newcommand{\T}{\mathbb{T}}
\newcommand{\Hh}{\mathcal{H}}
\newcommand{\iso}{\operatorname{Iso}}
\newcommand{\id}{\operatorname{id}}
\newcommand{\supp}{\operatorname{supp}}
\newcommand{\osupp}{\operatorname{supp}^\circ}
\newcommand{\rn}[1]{\lVert #1\rVert_r}
\newcommand{\mn}[1]{\lVert #1\rVert_\infty}
\newcommand{\go}{G^{(0)}}
\newcommand{\Io}[1]{\operatorname{Iso}(#1)^\circ}
\newcommand{\Stab}{\operatorname{Stab}}
\newcommand{\ess}{\operatorname{ess}}
\newcommand{\XX}{\overline{X}}
\newcommand{\Ss}{\overline{\sigma}}
\newcommand{\Oo}{\mathcal{O}}
\newcommand{\Dd}{\mathcal{D}}
\newcommand{\Cc}{\mathcal{C}}
\newcommand{\Kk}{\mathcal{K}}
\newcommand{\R}{\mathcal{R}}
\newcommand{\G}{G}
\newcommand{\Aut}{\operatorname{Aut}}
\newcommand{\NN}{{\N_0}}
\newcommand{\lsp}{\operatorname{span}}
\newcommand{\clsp}{\overline{\lsp}}

\theoremstyle{plain}
\theorembodyfont{\itshape}
\newmdtheoremenv[linecolor=blue,backgroundcolor=blue!10]{theorem}{Theorem}[section]
\newmdtheoremenv[linecolor=blue,backgroundcolor=blue!10]{lemma}[theorem]{Lemma}
\newmdtheoremenv[linecolor=blue,backgroundcolor=blue!10]{corollary}[theorem]{Corollary}
\newmdtheoremenv[linecolor=blue,backgroundcolor=blue!10]{proposition}[theorem]{Proposition}

\theorembodyfont{\normalfont}
\newmdtheoremenv[linecolor=blue,backgroundcolor=blue!10]{definition}[theorem]{Definition}
\newmdtheoremenv[linecolor=green,backgroundcolor=green!10]{example}[theorem]{Example}

\theoremstyle{break}
\newmdtheoremenv[linecolor=green,backgroundcolor=green!10]{examples}[theorem]{Examples}

\date{\today}
\author{Toke Meier Carlsen}
\title{$C^*$-rigidity of dynamical systems and étale groupoids}

\begin{document}
\maketitle\newpage
\setcounter{tocdepth}{2}
\tableofcontents

\setcounter{secnumdepth}{-1}
\chapter{Foreword}
These notes were written in connection with the masterclass \emph{Rigidity of $C^*$-algebras associated to dynamics} held at the University of Copenhagen October 16-20, 2017.

I wish to thank the organisers Sara E. Arklint, Kevin Aguyar Brix, and Søren Eilers for given me the oppertunity to talk about the exciting topic of $C^*$-rigidity of dynamical systems and étale groupoids.

There are doubtless errors in these notes. I would be grateful if you let me know if you find any typos, other errors, or things in the notes you think could be improved. 

\setcounter{secnumdepth}{2}
\chapter{Introduction}

$C^*$-rigidity of dynamical systems is the principle that dynamical systems can be recovered, up to a suitable notion of equivalence, from $C^*$-algebraic data associated to them.

Examples of this principle include characterisations of flip conjugacy and strong orbit equivalence of Cantor minimal systems in terms of crossed product $C^*$-algebras (see \cite{GiPuSk1995}), and of conjugacy, continuous orbit equivalence, and flow equivalence of topological Markov shifts in terms of Cuntz--Krieger algebras (see below).

I will in these notes explain how several of the above mentioned results, and also some new results, follow from a theorem about $C^*$-rigidity of étale groupoids, but let us first take a closer look at some of the $C^*$-rigidity results that will follow from this theorem.

\section{Homeormopshisms of compact Hausdorff spaces}

Let $X$ and $Y$ be compact Hausdorff spaces and let $\sigma : X \to X$ and $\tau : Y \to Y$ be homeomorphisms. Then $\sigma$ and $\tau$ are \emph{conjugate} if there is a homeomorphism $h:X\to Y$ such that $\tau\circ h=h\circ\sigma$, and \emph{flip-conjugate} if either $\sigma$ and $\tau$ are conjugate or $\sigma$ and $\tau^{-1}$ are conjugate. 

It is straightforward to show that if  $\sigma : X \to X$ and $\tau : Y \to Y$ are flip-conjugate, then there is a $*$-isomorphism between the corresponding crossed product $C^*$-algebras $C(X)\rtimes_\sigma\Z$ and $C(Y)\rtimes_\tau\Z$ that maps $C(X)$ onto $C(Y)$. By using a result in \cite{Boyle1983}, Giordano, Putnam, and Skau showed in their paper \cite{GiPuSk1995} that the converse holds for Cantor minimal systems (i.e., homeomorphisms $\sigma:X\to X$ where $X$ is a Cantor set such that $\{\sigma^n(x):n\in\Z\}$ is dense in $X$ for all $x\in X$). This result was generalised to the case where $X$ and $Y$ are compact metric spaces and  $\sigma$ and $\tau$ are topologically transitive (i.e., $\{\sigma^n(x):n\in\Z\}$ is dense in $X$ for some $x\in X$, and $\{\tau^n(y):n\in\Z\}$ is dense in $Y$ for some $y\in Y$) by Tomiyama in \cite{To1996}. 

These results were further generalised when Boyle and Tomiyama in their paper \cite{BoTo1998a} showed the following theorem for topologically free $\sigma$ and $\tau$ (i.e., when $\{x\in X:\sigma^n(x)=x\implies n=0\}$ is dense in $X$ and $\{y\in Y:\tau^n(y)=y\implies n=0\}$ is dense in $Y$), and then the general case was proven in \cite{CaRuSiTo2017}.

\begin{theorem}\label{thm:bt}
Suppose that $\sigma : X \to X$ and $\tau : Y \to Y$ are homeomorphisms of compact
Hausdorff spaces. Then there is a $*$-isomorphism from $C(X)\times_\sigma\Z$ to $C(Y)\times_\tau\Z$ that maps $C(X)$ onto $C(Y)$ if and only if there exist decompositions $X=X_1\sqcup X_2$ and $Y=Y_1\sqcup Y_2$ into disjoint open invariant sets such that $\sigma|X_1$ is conjugate to $\tau|Y_1$ and $\sigma|X_2$ is conjugate to $\tau^{-1}|Y_2$.

If $\sigma$ and $\tau$ are topologically transitive or $X$ and $Y$ are connected, then one of the sets $X_1$ and $X_2$ is empty and $\sigma$ and $\tau$ are flip-conjugate.
\end{theorem}

\section{Group actions and crossed products}

We shall now look at a result similar to Theorem~\ref{thm:bt} for topologically free actions of countable discrete groups on second-countable locally compact Hausdorff spaces.

Let $\Gamma$ be a countable discrete group acting continuously on the right of a locally compact second countable Hausdorff space $X$. We write $x\gamma$ for the action of $\gamma\in\Gamma$ on $x\in X$, and let $C_0(X)\rtimes_r\Gamma$ be the corresponding reduced crossed product $C^*$-algebra. We say that the action $\Gamma \curvearrowright X$ is \emph{topologically free} if the set $\{x\in X:x\gamma=x\implies \gamma=e\}$ (where $e$ is the identity element of $\Gamma$) is dense in $X$.

Two actions $\Gamma \curvearrowright X$ and $\Lambda \curvearrowright Y$ of countable discrete groups on second-countable locally compact Hausdorff spaces are \emph{continuously orbit equivalent} if there is a homeomorphism $h:X\to Y$ and continuous maps $\phi:X\times\Gamma\to\Lambda$ and $\eta:Y\times\Lambda\to\Gamma$ such that $h(x\gamma)=h(x)\phi(x,\gamma)$ for $x\in X$ and $\gamma\in\Gamma$, and $h^{-1}(y\lambda)=h^{-1}(y)\eta(y,\lambda)$ for $y\in Y$ and $\lambda\in\Lambda$.

Li gave in \cite{Li2016} the following characterisation of continuously orbit equivalence for topologically free actions of countable discrete groups on second-countable locally compact Hausdorff spaces.

\begin{theorem}\label{thm:group}
Let $\Gamma \curvearrowright X$ and $\Lambda \curvearrowright Y$ be topologically free actions of countable discrete groups on second-countable locally compact Hausdorff spaces. Then $\Gamma \curvearrowright X$ and $\Lambda \curvearrowright Y$ are continuously orbit equivalent if and only if there is a $*$-isomorphism from $C_0(X)\rtimes_r\Gamma$ to $C_0(Y)\rtimes_r\Lambda$ that maps $C_0(X)$ onto $C_0(Y)$.
\end{theorem}

Li also showed that by imposing further conditions of $X$, $Y$, $\Gamma$, $\Lambda$ and the actions, then it is possible to strengthen Theorem~\ref{thm:group}. He shows for example that if $X$ and $Y$ are connected, then there is a $*$-isomorphism from $C_0(X)\rtimes_r\Gamma$ to $C_0(Y)\rtimes_r\Lambda$ that maps $C_0(X)$ onto $C_0(Y)$ if and only if $\Gamma \curvearrowright X$ and $\Lambda \curvearrowright Y$ are \emph{conjugate}, i.e., if and only if there is a homeomorphism $h:X\to Y$ and a group isomorphism $\phi:\Gamma\to\Lambda$ such that $h(x\gamma)=h(x)\phi(\gamma)$ for $x\in X$ and $\gamma\in\Gamma$. Notice that if $\Gamma=\Z=\Lambda$, then $\Gamma \curvearrowright X$ and $\Lambda \curvearrowright Y$ are conjugate if and only if the homeomorphisms $x\mapsto x1$ and $y\mapsto y1$ are flip-conjugate.

We shall in Chapter~\ref{chapter:dyn} see how Theorem~\ref{thm:group} can be generalised to actions that are not necessarily topologically free.

\section{Cuntz--Krieger algebras and topological Markov shifts}

We end this chapter by looking at how conjugacy, eventually conjugacy, continuous orbit equivalence, and flow equivalence of topological Markov shifts can be characterised in terms of Cuntz--Krieger algebras.

\emph{Cuntz--Krieger algebras} is an important class of $C^*$-algebras introduced by Cuntz and Krieger in \cite{CuKr1980}, and further studied by Cuntz in \cite{Cu1981}. The definition of a Cuntz--Krieger algebra as a universal $C^*$-algebra is from \cite{HuRa1997}.

\begin{definition}
Let $A=(A_{ij})_{i,j=1}^n$ be an $n\times n$ matrix with entries in $\{0,1\}$ and with no zero rows and no zero columns. The Cuntz--Krieger algebra of $A$ is the universal $C^*$-algebra $\Oo_A$ generated by $n$ partial isometries $S_1,\dots,S_n$ satisfying
\begin{equation}\label{eq:ck}
S_i^*S_i=\sum_{i=1}A_{ij}S_jS_j^*\text{ and }\sum_{i=1}^nS_iS_i^*=1.
\end{equation}
\end{definition}

That $\Oo_A$ is the \emph{universal} $C^*$-algebra generated by $n$ partial isometries $S_1$,\linebreak\dots, $S_n$ satisfying \eqref{eq:ck} means that $\mathcal{O}_A$ is generated by partial isometries $S_1,\dots,S_n$ satisfying \eqref{eq:ck}, and that if $s_1,\dots,s_n$ are $n$ partial isometries in a unital $C^*$-algebra $\mathcal{A}$ satisfying $s_i^*s_i=\sum_{i=1}A_{ij}s_js_j^*$ and $\sum_{i=1}^ns_is_i^*=1$, then there is a $*$-homomorphism from $\Oo_A$ to $\mathcal{A}$ sending $S_i$ to $s_i$ for $i\in\{1,\dots,n\}$.

Already in \cite{CuKr1980}, Cuntz and Krieger described how the structure of the $C^*$-algebra $\Oo_A$ is closely related to the structure of a certain dynamical system $(\XX_A,\Ss_A)$ constructed from $A$.

\begin{definition}
Let $A=(A_{ij})_{i,j=1}^n$ be an $n\times n$ matrix with entries in $\{0,1\}$ and with no zero rows and no zero columns. Let 
\[\XX_A:=\left\{(x)_{k\in\Z}\in\{1,\dots,n\}^\Z:A_{x_kx_{k+1}}=1\text{ for all }k\in\Z\right\},\]
and equip $\XX_A$ with the subspace topology of $\{1,\dots,n\}^\Z$ where the latter is given the product topology of the discrete topology of $\{1,\dots,n\}$. 

Define $\Ss_A:\XX_A\to\XX_A$ by $\Ss_A((x)_{k\in\Z})=(y_k)_{k\in\Z}$ where $y_k=x_{k+1}$ for $k\in\Z$. 
\end{definition}

It is easy to check that $\XX_A$ is both open and closed in $\{1,\dots,n\}^\Z$. Since $\{1,\dots,n\}^\Z$ is a compact Hausdorff space, it follows that also $\XX_A$ is compact and Hausdorff. It is also zero-dimensional (i.e., its topology has a basis consisting of clopen sets) and metrizable (one can for instance use the metric given by $$d((x)_{k\in\Z},(y_k)_{k\in\Z})=2^{-\min\{|k|:x_k\ne y_k\}}$$ if $(x)_{k\in\Z}\ne (y_k)_{k\in\Z}$), but it is usually easier just to use the topology instead of a metric.

It is straightforward to check that $\Ss_A$ is a homeomorphism from $\XX_A$ to itself. The dynamical $(\XX_A,\Ss_A)$ is sometimes called a \emph{topological Markov shift}. It is an example of a \emph{shift of finite type}. Actually, any shift of finite type is conjugate (two dynamical systems $(X,\sigma)$ and $(Y,\tau)$ are \emph{conjugate} if there is a homeomorphism $\phi:X\to Y$ such that $\phi\circ\sigma=\tau\circ\phi$) to $(\XX_A,\Ss_A)$ for some square matrix $A$ with entries in $\{0,1\}$  and with no zero rows and no zero columns (see for example \cite{LiMa1995}).

We have already mentioned what it means for two dynamical systems to be conjugate. For topological Markov shifts there is another, weaker, form of equivalence we will mention here, and that is \emph{flow equivalence}. We will not define flow equivalence here (see for example \cite{LiMa1995}), but instead use the following characterisation of flow equivalence of zero-dimensional dynamical system given by Parry and Sullivan in \cite{PaSu1975} (see also \cite{BoCaEi2017}).


\begin{definition}
Let $X$ be a compact Hausdorff space, and $\sigma:X\to X$ a homeomorphism. A \emph{discrete cross section} to $(X,\sigma)$ is a pair $(X',\sigma')$ where $X'$ is a closed subset of $X$ and $\sigma'$ is a homeomorphism from $X'$ to itself such that there is a continuous function $r:X'\to\N$ such that $\sigma'(x)=\sigma^{r(x)}(x)$ and $\sigma^{j}(x)\notin X'$ for all $x\in X'$ and all $j\in\{1,2,\dots,r(x)-1\}$, and $X=\{\sigma^n(x):x\in X',\ n\in\N\}$.
\end{definition}

The map $r:X'\to\N$ above is often called the \emph{return time map} of $X'$, and $r(x)$ is called the \emph{return time} of $x$.

\begin{theorem}[Parry and Sullivan]
Let $X$ and $Y$ be compact metrizable zero-dimensional spaces and let $\sigma:X\to X$ and $\tau:Y\to Y$ be homeomorphisms. Then $(X,\sigma)$ and $(Y,\tau)$ are flow equivalent if and only if there is a discrete cross section $(X',\sigma')$ of $(X,\sigma)$ and a discrete cross section $(Y',\tau')$ of $(Y,\tau)$ such that $(X',\sigma')$ and $(Y',\tau')$ are conjugate.
\end{theorem}

The topological Markov shift $(\XX_A,\Ss_A)$ is an example of a two-sided subshift. From such a two-sided shift space, one can construct a one-sided subshift (see \cite{Ki1998}).

\begin{definition}
Let $A=(A_{ij})_{i,j=1}^n$ be an $n\times n$ matrix with entries in $\{0,1\}$ and with no zero rows and no zero columns. Let 
\[X_A:=\left\{(x)_{k\in\N_0}\in\{1,\dots,n\}^{\N_0}:A_{x_kx_{k+1}}=1\text{ for all }k\in\N_0\right\},\]
and equip $X_A$ with the subspace topology of $\{1,\dots,n\}^{\N_0}$ where the latter is given the product topology of the discrete topology of $\{1,\dots,n\}$. 

Define $\sigma_A:X_A\to X_A$ by $\sigma_A((x)_{k\in\N_0})=(y_k)_{k\in\N_0}$ where $y_k=x_{k+1}$ for $k\in\N_0$. 
\end{definition}

As $\XX_A$, the space $X_A$ is compact, metrizable and zero-dimensional. An important difference between $(\XX_A,\Ss_A)$ and $(X_A,\sigma_A)$ is that $\sigma_A$ is not necessarily a homeomorphism, but it is a surjective local homeomorphism (i.e., there is for each $x\in X_A$ an open neighbourhood $U$ of $x$ such that $\sigma_A(U)$ is open and the restriction $\sigma_A|_U$ of $\sigma_A$ to $U$ is a homeomorphism from $U$ to $\sigma_A(U)$).

As with two-sided topological Markov shifts, we say that two one-sided topological Markov shifts $(X_A,\sigma_A)$ and $(X_B,\sigma_B)$ are conjugate if there is a homeomorphism $h:X_A\to X_B$ such that $h\circ\sigma_A=\sigma_B\circ h$. We shall look at two other, and weaker, forms of equivalence of one-sided subshifts.

\begin{definition}
Let $A$ and $B$ be two square matrices with entries in $\{0,1\}$ and with no zero rows and no zero columns. 

We say that $(X_A,\sigma_A)$ and $(X_B,\sigma_B)$ are \emph{eventually conjugate} if there is a homeomorphism $h:X_A\to X_B$ and an $k\in\N_0$ such that $$\sigma_B^{k+1}(h(x))=\sigma_B^{k}(h(\sigma_A(x)))$$ for $x\in X_A$, and $$\sigma_A^{k+1}(h^{-1}(y))=\sigma_A^{k}(h^{-1}(\sigma_B(y)))$$ for $y\in X_B$.

We say that $(X_A,\sigma_A)$ and $(X_B,\sigma_B)$ are \emph{continuously orbit equivalent} if there is a homeomorphism $h:X_A\to X_B$ and an continuous functions $k_A,l_A:X_A\to\N_0$ and $k_B,l_B:X_B\to\N_0$ such that $$\sigma_B^{l_A(x)}(h(x))=\sigma_B^{k_A(x)}(h(\sigma_A(x)))$$ for $x\in X_A$, and $$\sigma_A^{l_B(y)}(h^{-1}(y))=\sigma_A^{k_B(y)}(h^{-1}(\sigma_B(y)))$$ for $y\in X_B$.
\end{definition}

Notice that conjugacy implies eventually conjugacy, and that eventually conjugacy implies continuous orbit equivalence.

Let $A=(A_{ij})_{i,j=1}^n$ be an $n\times n$ matrix with entries in $\{0,1\}$ and with no zero rows and no zero columns. As in \cite{CuKr1980}, we denote by $\Dd_A$ the $C^*$-subalgebra of $\Oo_A$ generated by elements of the form $S_{i_0}S_{i_1}\dots S_{i_j}S_{i_j}^*\dots S_{i_1}^*S_{i_0}^*$ where $j\in\N_0$ and $i_0,i_1,\dots,i_j\in\{1,\dots,n\}$. We then have that $\Dd_A$ is isomorphic to $C(X_A)$ by an isomorphism $\omega$ that maps 
\[S_{i_0}S_{i_1}\dots S_{i_j}S_{i_j}^*\dots S_{i_1}^*S_{i_0}^*\] 
to the indicator function of the set $\{(x_k)_{k\in\N_0}\in X_A:x_0=i_0,\dots,x_j=i_j\}$. For $t\in\T$, let $\lambda_t^A$ be the automorphism of $\Oo_A$ given by $\lambda_t^A(S_i)=tS_i$ for $i\in\{1,\dots,n\}$, and let $\tau_A:\Oo_A\to\Oo_A$ be the completely positive map defined by $\tau_A(X)=\sum_{i,j=1}^nS_iXS_j^*$. Let $\Kk$ denote the $C^*$-algebra of compact operators on a separable infinite-dimensional Hilbert space, and let $\Cc$ denote a maximal abelian $C^*$-subalgebra of $\Kk$.

We are now ready to state the above mentioned characterisation of conjugacy, eventually conjugacy, continuous orbit equivalence and flow equivalence of topological Markov shifts in terms of Cuntz--Krieger algebras.

\begin{theorem}\label{thm:ck}
Let $A$ and $B$ be two square matrices with entries in $\{0,1\}$ and with no zero rows and no zero columns. 
\begin{enumerate}
	\item $(X_A,\sigma_A)$ and $(X_B,\sigma_B)$ are conjugate if and only if there is a $*$-isomor\-phism $\psi:\Oo_A\to\Oo_B$ such that $\psi(\Dd_A)=\Dd_B$ and $\psi\circ\tau_A=\tau_B\circ\psi$.
	\item $(X_A,\sigma_A)$ and $(X_B,\sigma_B)$ are eventually conjugate if and only if there is a $*$-isomorphism $\psi:\Oo_A\to\Oo_B$ such that $\psi(\Dd_A)=\Dd_B$ and $\psi\circ\lambda_t^A=\lambda_t^B\circ\psi$ for all $t\in\T$.
	\item $(X_A,\sigma_A)$ and $(X_B,\sigma_B)$ are continuously orbit equivalent if and only if there is a $*$-isomorphism $\psi:\Oo_A\to\Oo_B$ such that $\psi(\Dd_A)=\Dd_B$.
	\item $(\XX_A,\Ss_A)$ and $(\XX_B,\Ss_B)$ are conjugate if and only if there is a $*$-isomor\-phism $\psi:\Oo_A\otimes\Kk\to\Oo_B\otimes\Kk$ such that $\psi(\Dd_A\otimes\Cc)=\Dd_B\otimes\Cc$ and $\psi\circ (\lambda_t^A\otimes\id)=(\lambda_t^B\otimes\id)\circ\psi$ for all $t\in\T$.
	\item $(\XX_A,\Ss_A)$ and $(\XX_B,\Ss_B)$ are flow equivalent if and only if there is a $*$-isomorphism $\psi:\Oo_A\otimes\Kk\to\Oo_B\otimes\Kk$ such that $\psi(\Dd_A\otimes\Cc)=\Dd_B\otimes\Cc$.
\end{enumerate}
\end{theorem}

Theorem~\ref{thm:ck}(1) was proven in \cite[Theorem 3.3]{CarlsenBrix}.\linebreak Theorem~\ref{thm:ck}(2) was proven for irreducible matrices satisfying condition (I) of \cite{CuKr1980} in \cite{Ma2016c}, and in general in the paper \cite{CaRo2017}. It is related to \cite[Proposition 2.17]{CuKr1980}. 

Theorem~\ref{thm:ck}(3) was proven for irreducible matrices satisfying condition (I) of \cite{CuKr1980} in \cite{Ma2010}. The general case follows from \cite[Corollary 4.6]{CaWi2016} and \cite[Theorem 5.3]{AREIRU2017} (or \cite[Theorem 5.1]{BrCaWh2017} and \cite[Theorem 4.2]{CaWi2016}). 

Theorem~\ref{thm:ck}(4) was proven in \cite[Corollary 5.2]{CaRo2017} (the ``only if'' part of it was proven for irreducible matrices satisfying condition (I) of \cite{CuKr1980} in \cite[Theorem 3.8]{CuKr1980}, and for matrices with the property that both themselves and their transposes satisfying condition (I) of \cite{CuKr1980} in \cite[Theorem 2.3]{Cu1981}). 

Theorem~\ref{thm:ck}(5) was proven for irreducible matrices satisfying condition (I) of \cite{CuKr1980} in \cite{MaMa2016} and in general in \cite{CaEiOr2016} (the ``only if'' part of it was proven for irreducible matrices satisfying condition (I) of \cite{CuKr1980} in \cite[Theorem 4.1]{CuKr1980}, and for matrices with the property that both themselves and their transposes satisfying condition (I) of \cite{CuKr1980} in the paper \cite[Theorem 2.4]{Cu1981}).

We shall in Chapter~\ref{chapter:dyn} see how Theorem~\ref{thm:ck} can be proved by using grou\-poids.

\chapter{Étale groupoids and their $C^*$-algebras}\label{chapter:groupoids}

$C^*$-algebras of topological groupoids were introduced in \cite{Re1980a} and have since then been intensively studied. Many interesting classes of $C^*$-algebras can be constructed as $C^*$-algebras of topological groupoids, which therefore provides a unified framework for studying things like the ideal structure, K-theory, KMS-states, pure infiniteness, isomorphisms, and Morita equivalence for these classes of $C^*$-algebras. 

I will in this chapter give a short introduction to étale groupoids and the reduced $C^*$-algebra of an étale Hausdorff groupoid. The subclass of étale groupoids plays in many ways the same role in the class of topological groupoids as the subclass of discrete groups plays in the class of topological groups. 

I recommend \cite{Re1980a}, \cite{Pa1999}, and \cite{Sims2017} for more on groupoids and their $C^*$-algebras.

\section{Groupoids}

\begin{definition}\label{def:groupoid}
A \emph{groupoid} is a small category\footnote{That a category is \emph{small} means that both its class of objects and its class of morphisms are sets.} in which every morphism has an inverse.
\end{definition}
If $G$ is a groupoid, then we write $G^{(0)}$ (or sometimes just $G^0$) to denote its set of objects and $G^{(1)}$ to denotes its set of morphisms. Usually, at least when using groupoids in connection with $C^*$-algebras, one identify an object with its corresponding identity morphism. By doing this, $G^{(0)}$ is then regarded as a subset of $G^{(1)}$, and one then often just write $G$ instead of $G^{(1)}$. We shall use this convention throught these notes.

It is also customary to call the domain of a morphism $\eta$ for its \emph{source} and denote it by $s(\eta)$, and to call the codomain of $\eta$ for its \emph{range} and denote it by $r(\eta)$. We then have that $s$ and $r$ are maps from $G$ to $G^{(0)}$. The composition $\eta_1\eta_2$ of two elements $\eta_1$ and $\eta_2$ in $G$ is then defined if and only if $s(\eta_1)=r(\eta_2)$. Often the notation $G^{(2)}:=\{(\eta_1,\eta_2)\in G\times G:s(\eta_1)=r(\eta_2)\}$ is used for the set of composable pairs. It is also customary to call the composition $\eta_1\eta_2$ for the \emph{product} of $\eta_1$ and $\eta_2$.

The inverse of an element $\eta\in G$ is denoted by $\eta^{-1}$. We thus have that $\eta\mapsto\eta^{-1}$ is a map from $G$ to $G$.

The sets $G$, $G^{(0)}$, and $G^{(2)}$, the range and source maps $r$ and $s$, and the product and the inverse map have the following properties.
\begin{enumerate}[\indent (1)]
	\item $r(x)=x=s(x)$ for all $x\in G^{(0)}$.
	\item $r(\eta)\eta=\eta=\eta s(\eta)$ for all $\eta\in G$.
	\item $r(\eta^{-1})=s(\eta)$ and $s(\eta^{-1})=r(\eta)$ for all $\eta\in G$.
	\item $\eta^{-1}\eta=s(\eta)$ and $\eta\eta^{-1}=r(\eta)$ for all $\eta\in G$.
	\item $r(\eta_1\eta_2)=r(\eta_1)$ and $s(\eta_1\eta_2)=s(\eta_2)$ for all $(\eta_1,\eta_2)\in G^{(2)}$.
	\item $(\eta_1\eta_2)\eta_3=\eta_1(\eta_2\eta_3)$ whenever $(\eta_1,\eta_2), (\eta_2,\eta_3)\in G^{(2)}$.
\end{enumerate}
Conversely; if $G$ is a set, $G^{(0)}$ is a subset of $G$, $r$ and $s$ are maps from $G$ to $G^{(0)}$, $G^{(2)}:=\{(\eta_1,\eta_2)\in G\times G:s(\eta_1)=r(\eta_2)\}$, and $(\eta_1,\eta_2)\mapsto \eta_1\eta_2$ is a map from $G^{(2)}$ to $G$, and $\eta\mapsto\eta^{-1}$ is a map from $G$ to $G$ such that (1)--(6) hold; then there is a groupoid such that $G$ is its set of morphism, $G^{(0)}$ is its set of objects, $s(\eta)$ is the domain, $r(\eta)$ is the codomain, and $\eta^{-1}$ is the inverse of an morphisms $\eta$, and the composition of two morphisms $\eta_1$ and $\eta_2$ for which $s(\eta_1)=r(\eta_2)$ is $\eta_1\eta_2$. 
Thus, instead of Definition~\ref{def:groupoid}, an alternative, but equivalent way of defining a groupoid is to say that a groupoid consists of a set $G$, a subset $G^{(0)}$ of $G$, maps $r$ and $s$ from $G$ to $G^{(0)}$, a map $(\eta_1,\eta_2)\mapsto \eta_1\eta_2$ from $\{(\eta_1,\eta_2)\in G\times G:s(\eta_1)=r(\eta_2)\}$ to $G$, and a map $\eta\mapsto\eta^{-1}$ from $G$ to $G$ such that (1)--(6) hold. 

There is an alternative way of characterising a groupoid where one specify the set $G$, the set $G^{(2)}$, the product, and the inverse, and then define the set $G^{(0)}$ and the maps $r$ and $s$ from this (see \cite[Definition 1.1]{Ha1978}).

\begin{example}\label{it:group} 
Let $G$ be a group and let $e$ be its identity. Then $G$ is a groupoid with $G^{(0)}:=\{e\}$, and the product and inverse given by the group operations.
\end{example}

\begin{example}\label{it:space} 
Let $X$ be a set. Then $X$ is a groupoid with $X^{(0)}:=X$, $r$ and $s$ the identity maps, the product defined by $(x,x)\mapsto x$, and the inverse defined by $x^{-1}=x$.
\end{example}

\begin{example}\label{it:bundle} 
Let $(E,X,\pi)$ be a group bundle, i.e., $E$ and $X$ are sets, $\pi$ is a surjevtive map from $E$ to $X$, and $\pi^{-1}(x)$ is a group for each $x\in X$. Then $E$ is a groupoid with $E^{(0)}=\{e_x:x\in X\}$, where for each $x\in X$, $e_x$ is the identity of $\pi^{-1}(x)$; $r(\eta)=s(\eta)=e_{\pi(\eta)}$ and $\eta^{-1}$ is the inverse of $\eta$ in $\pi^{-1}(\pi(\eta))$; and the product of $\eta_1$ and $\eta_2$ is the product of $\eta_1$ and $\eta_2$ in $\pi^{-1}(\pi(\eta_1))=\pi^{-1}(\pi(\eta_2))$.
\end{example}

\begin{example} 
Let $X$ be a set and $\sim$ an equivalence relation on $X$. Let $G:=\{(x,y)\in X\times X:x\sim y\}$, let $\go:=\{(x,x)\in G:x\in X\}$ which we identify with $X$, and define $r,s:G\to X$ by $r(x,y)=x$ and $s(x,y)=y$. For $(x_1,y_1), (x_2,y_2)$ with $y_2=x_1$, let $(x_1,y_1)(x_2,y_2)=(x_1,y_2)$; and let $(x,y)^{-1}=(y,x)$ for $(x,y)\in G$. Then $G$ is a groupoid.
\end{example}

Let $G$ be a groupoid. For $x\in G^{(0)}$ we let $xG:=\{\eta\in G:r(\eta)=x\}$, $Gx:=\{\eta\in G:s(\eta)=x\}$, and $xGx:=xG\cap Gx=\{\eta\in G:s(\eta)=r(\eta)=x\}$ (often, the notation $G^x$ is used instead of $xG$, the notation $G_x$ is used instead of $Gx$, and the notation $G_x^x$ is used instead of $xGx$). We also let $\iso(G):=\bigcup_{x\in G^{(0)}}xGx=\{\eta\in G:s(\eta)=r(\eta)\}$. Notice that $xGx$ and $\iso(G)$ are subgroupoids of $G$, and that $xGX$ is the groupoid corresponding to a group (cf. Example~\ref{it:group}, and that $\iso(G)$ is the groupoid corresponding to a group bundle (cf. Example~\ref{it:bundle}). The subgroupoid $xGx$ is called the \emph{isotropy of $x$} (or the \emph{isotropy group of $x$}), and the subgroupoid $\iso(G)$ is called the \emph{isotropy of $G$} (or the \emph{isotropy bundle of $G$}).

The \emph{orbit} of an $x\in G^{(0)}$ is the set $\{r(\eta):\eta\in Gx\}$. In order words, $x\in G^{(0)}$ and $x'\in G^{(0)}$ belongs to the same orbit, if and only if there is an $\eta\in G$ such that $r(\eta)=x$ and $s(\eta)=x'$. Notice that we in that case have that $\eta'\mapsto \eta\eta'\eta^{-1}$ is an isomorphism from $x'Gx'$ to $xGx$.

A subset $A$ of a groupoid $G$ is called a \emph{bisection} if the restrictions of $r$ and $s$ to $A$ are both injective.

\section{Étale groupoids}

\begin{definition}
A \emph{topological groupoid} is a groupoid $G$ endowed with a topology under which the maps $r$ are $s$ are continous maps from $G$ to $G^{(0)}$ with respect to the relative topology on $G^{(0)}$, the map $\eta\mapsto\eta^{-1}$ is a continous map from $G$ to $G$, and the map $(\eta_1,\eta_2)\mapsto \eta_1\eta_2$ is a continous map from $G^{(2)}$ to $G$ with respect to the relative topology of the product topology on $G^{(2)}\subseteq G\times G$.
\end{definition}

In these notes, the topological groupoids we consider are always Hausdorff, but there are plenty of natural occuring examples of topological groupoids that are not Hausdorff. If a topological groupoid $G$ is Hausdorff, then $G^{(0)}$ is a closed subset of $G$.

\begin{definition}
An \emph{étale} groupoid is a locally compact topological groupoid $G$ such that $G^{(0)}$ is an open subset of $G$, $G^{(0)}$ is Hausdorff in the relative topology, and $r:G\to G^{(0)}$ (equivalently $s:G\to G^{(0)}$) is a local homeomorphism (i.e., each $\eta\in G$ has an open neighbourhood $U$ such that $r(U)$ is open and the restriction of $r$ to $U$ is a homeomorphism from $U$ to $r(U)$).
\end{definition}

\begin{example}
If $G$ is any groupoid, then $G$ becomes an étale groupoid if we equip it with the discrete topology.
\end{example}

\begin{example}
Let $\Gamma$ be a topological group acting continuously on the right on a topological space $X$. We write $x\gamma$ for the action of $\gamma$ on $x$. Let 
\[
X \rtimes \Gamma := X\times\Gamma
\]
and equip $X \rtimes \Gamma$ with the product topology. Let $(X \rtimes \Gamma)^{(0)}:=X \times \{e\}$, which we identify with $X$, and define $r,s:X \rtimes \Gamma\to X$ by $r(x,\gamma) = x$ and $s(x,\gamma) = x\gamma$. Then $\big((x_1, \gamma_1), (x_2, \gamma_2)\big) \in (X \rtimes \Gamma)^{(2)}$ if and only if $\gamma_2 = x_1\gamma_1$, in which case we let $(x_1,\gamma_1)(x_1\gamma_1, \gamma_2) := (x_1, \gamma_1\gamma_2)$. We also let $(x, \gamma)^{-1} := (x\gamma, \gamma^{-1})$.  

Then is a $X \rtimes \Gamma$ is topological groupoid, which is étale if and only if $\Gamma$ is discrete and $X$ is locally compact and Hausdorff.
\end{example}

\begin{example}
Let $X$ be a locally compact Hausdorff space, let $U$ and $V$ be open subsets of $X$ and let $\sigma:U\to V$ be a surjective local homeomorphism (i.e., there is for each $x\in U$ an open neighbourhood $O\subseteq U$ of $x$ such that $\sigma(O)$ is open and the restriction of $\sigma$ to $O$ is a homeomorphism from $O$ to $\sigma(O)$). We inductively define $U_n$, $V_n$, and $\sigma^n:U_n\to V_n$ in the following way. Let $U_0:=V_0:=X$ and let $\sigma^0:X\to X$ be the identity map. Let $U^1:=U$, $V^1:=V$ and $\sigma^1:=\sigma$. For $n>1$, let $U^n:=\sigma^{-1}(U_{n-1}\cap V)$, $V_n=\sigma(V_{n-1}\cap U)$, and $\sigma^n:=\sigma^{n-1}\circ\sigma|_{U_n}$. Then $U_n$ and $V_n$ are open and $\sigma^n:U_n\to V_n$ is a local homeomorphism for each $n\in\mathbb{N}_0$.

Let 
\[G(X,\sigma):=\bigcup_{n,m\in\mathbb{N}_0}\{(x,m-n,y)\in U_m\times\{m-n\}\times U_n: \sigma^m(x)=\sigma^n(y)\}\]
let $G(X,\sigma)^{(0)}:=\{(x,0,x)\in G(X,\sigma):x\in X\}$ which we identify with $X$, and define $r,s:G(X,\sigma)\to X$ by $r(x,k,y)=x$ and $s(x,k,y)=y$. For $(x_1,k_1,y_1)$, $(x_2,k_2,y_2)\in G(X,\sigma)$ with $y_1=x_2$, let $(x_1,k_1,y_1)(x_2,k_2,y_2)=(x_1,k_1+k_2,y_2)\in G(X,\sigma)$; and for $(x,k,y)\in G(X,\sigma)$, let $(x,k,y)^{-1}:=(y,-k,x)\in G(X,\sigma)$. Then $G(X,\sigma)$ is a groupoid.

For $m,n\in\mathbb{N}_0$ and open subsets $A$ of $U_m$ and $B$ of $U_n$ for which $\sigma^m|_A:A\to\sigma^m(A)$ and $\sigma^n|_B:B\to\sigma^n(B)$ are homeomorphisms, let
\[Z(A,m,n,B):=\{(x,m-n,y)\in G(X,\sigma): x\in A,\ y\in B,\ \sigma^m(x)=\sigma^n(y)\}.\]
Then the collection 
\begin{align*}
\{Z(A,m,n,B):&m,n\in\mathbb{N}_0,\ 
A\text{ is an open subset of }U_m,\\ 
&B\text{ is an open subset of }U_n,
\\&\sigma^m|_A:A\to\sigma^m(A)\text{ and }\sigma^n|_B:B\to\sigma^n(B)\text{ are homeomorphisms}\}
\end{align*}
is a basis for a topology on $G(X,\sigma)$ that makes $G(X,\sigma)$ an étale groupoid.
\end{example}

If $G$ is en étale groupoid, then is has a basis consisting of open bisections, and $xG$ and $Gx$ are, for each $x\in G^{(0)}$, discrete in the relative topologies.

\section{The $C^*$-algebra of an étale groupoid}
Let $G$ be a locally compact Hausdorff étale groupoid. If $f,g\in C_c(G)$ and $\eta\in G$, then the set $\{(\eta_1,\eta_2)\in G^{(2)}:\eta_1\eta_2=\eta,\ f(\eta_1)g(\eta_2)\ne 0\}$ is finite. We can therefore define a function $f\ast g:G\to\C$ by 
$$f\ast g(\eta):=\sum_{\eta_1\eta_2=\eta}f(\eta_1)g(\eta_2).$$ 
It is not difficult to check that $f\ast g\in C_c(G)$. The complex vector space $C_c(G)$ is a $*$-algebra with multiplication given by $\ast$ and involution given by $f^*(\eta)=\overline{f(\eta^{-1})}$.

As with group algebras, there are two pre-$C^*$-norms on $C_c(G)$; the universal norm and the reduced norm. By completing $C_c(G)$ with respect to the universal norm, we obtains the full $C^*$-algebra $C^*(G)$ of $G$, and by completing $C_c^(G)$ with respect to the reduced norm, we obtains the reduced $C^*$-algebra $C_r^*(G)$ of $G$. For many groupoids, the universal and the reduced norms coincide. In these notes, we will not consider the universal norm or the full $C^*$-algebra $C^*(G)$ of $G$. The reduced norm on $C_c(G)$ is defined in the following way. For each $x\in G^{(0)}$, there is a $*$-representation of $\pi_x:C_c(G)\to \mathcal{B}(l^2(Gx))$ such that $\pi_x(f)\delta_\eta=\sum_{\eta'\in Gr(\eta)}f(\eta')\delta_{\eta'\eta}$. The reduced $C^*$-algebra $C_r^*(G)$ of $G$ is the completion of $\bigoplus_{x\in G^{(0)}}\pi_x(C_c(G))$ in $\bigoplus_{x\in G^{(0)}}\mathcal{B}(l^2(Gx))$. The norm of $C^*_r(G)$ is called the \emph{reduced norm} and is denoted by $\rn{\cdot}$. If $G$ is second-countable, then $C^*_r(G)$ is separable.

The map $f\mapsto \bigoplus_{x\in G^{(0)}}\pi_x(f)$ from $C_c(G)$ to $C^*_r(G)$ is injective, and $\mn{f}\le \rn{\bigoplus_{x\in G^{(0)}}\pi_x(f)}$ for all $f\in C_c(G)$. We will throughout these notes identify $C_c(G)$ with its image under this map, and thus write $f$ instead of $\bigoplus_{x\in G^{(0)}}\pi_x(f)$ and regard $C_c(G)$ as a subalgebra of $C^*_r(G)$. If $f\in C_c(G)$ is supported on an open bisection, then $\mn{f}=\rn{f}$. 

There is an injective map $j:C^*_r(G)\to C_0(G)$ with the property that $j(f)(\eta)=\langle \pi_{s(\eta)}(f)\delta_{s(\eta)},\delta_\eta\rangle$ for $f\in C^*_r(G)$ and $\eta\in G$, where $\langle\cdot,\cdot\rangle$ is the inner product on $l^2(Gx)$. Moreover, $\mn{j(f)}\le \rn{f}$. If $f\in C_c(G)$, then $j(f)=f$.

We let $\iso(G)^\circ$ denote the interior of $\iso(G)$. Then $\iso(G)^\circ$ is a subgroupoid of $G$ with unit space $G^{(0)}$. 

If $G$ is a locally compact Hausdorff étale groupoid and $H$ is an open subgroupoid of $G$ (i.e. $H$ is an open subset of $G$ such that $\eta\in H\implies \eta^{-1}\in H$ and $(\eta_1,\eta_2)\in G^{(2)}\cap H\times H\implies \eta_1\eta_2\in H$), then $H$ is a locally compact Hausdorff étale groupoid with the operation it inherits from $G$. The inclusion of $C_c(H)$ into $C_c(G)$ extends to an inclusion of $C^*_r(H)$ into $C^*_r(G)$. Since $G^{(0)}$ and $\iso(G)^\circ$ are open subgroupoids of $G$, we will regard $C_0(G^{(0)})=C^*_r(G^{(0)})$ and $C^*_r(\iso(G)^\circ)$ to be $C^*$-subalgebras of $C^*_r(G)$. We then have that $C_0(G^{(0)})\subseteq C^*_r(\iso(G)^\circ)$. Morover, since the unit spaces of $G^{(0)}$, $\iso(G)^\circ$, and $G$ are equal, it follows that $C_0(G^{(0)})$ contains an approximate unit of $C^*_r(G)$ (and therefore also of $C^*_r(\iso(G))$).

If $\iso(G)^\circ$ is abelian (i.e., $\eta_1\eta_2=\eta_1\eta_2$ if $s(\eta_1)=r(\eta_1)=s(\eta_2)=r(\eta_2)$ and $\eta_1,\eta_2\in \iso(G)^\circ$), then $C^*_r(\iso(G)^\circ)$ is equal to the relative commutant $\{f\in C^*_r(G): fg=gf\text{ for all }g\in C_0(G^{(0)})\}$ of $C_0(G^{(0)})$ in $C^*_r(G)$. This means that $C^*_r(\iso(G)^\circ)$ is a $C_0(G^{(0)})$-algebra (i.e., there is an inclusion of $C_0(G^{(0)})$ into the center of $C^*_r(\iso(G)^\circ)$ which is nondegenerate in the sense that $\operatorname{span}\{gf:g\in C_0(G^{(0)}),\ f\in C^*_r(\iso(G)^\circ)\}$ is dense in $C^*_r(\iso(G)^\circ)$). We then have for each $x\in G^{(0)}$ an ideal $J_x:=\overline{\operatorname{span}}\{gf:g\in C_0(G^{(0)}),\ g(x)=0,\ f\in C^*_r(\iso(G)^\circ)\}$. We denote by $\pi_x$ the corresponding quotient map from $C^*_r(\iso(G)^\circ)$ to $C^*_r(\iso(G)^\circ)/J_x$. Then the map $x\mapsto \lVert\pi_x(f)\rVert$ is continuous for each $f\in C^*_r(\iso(G)^\circ)$, and there is for each $x\in G^{(0)}$ an isomorphism between $\pi_x(C^*_r(\iso(G)^\circ))$ and $C^*(xGx^\circ)$ that maps $\pi_x(f)$ to $f(x)1_{C^*(xGx^\circ)}$ for each $f\in C^*_r(\iso(G)^\circ)$.

\section{Cocycles and coactions}
Let $\Gamma$ be a discrete group. A \emph{cocycle} from $G$ to $\Gamma$ is a map $c:G\to\Gamma$ such that $c(\eta^{-1})=c(\eta)^{-1}$ for $\eta\in G$, and $c(\eta_1\eta_2)=c(\eta_1)c(\eta_2)$ for $(\eta_1,\eta_2)\in G^{(2)}$. A cocycle $c:G\to\Gamma$ induces a $\Gamma$-\emph{grading} $\{c^{-1}(\gamma)\}_{\gamma\in\Gamma}$ of $G$ (i.e., $\bigcup_{\gamma\in\Gamma}c^{-1}(\gamma)=G$, $c^{-1}(\gamma_1)\cap c^{-1}(\gamma_2)=\emptyset$ for $\gamma_1\ne\gamma_2$, and $\eta_1\eta_2\in c^{-1}(\gamma_1\gamma_2)$ if $(\eta_1,\eta_2)\in G^{(2)}$, $\eta_1\in c^{-1}(\gamma_1)$, and $\eta_2\in c^{-1}(\gamma_2)$). It also induces a reduced coaction $\delta_c:C^*_r(G)\to C^*_r(G)\otimes C^*_r(\Gamma)$ such that $\delta_c(f)=f\otimes \lambda_\gamma$ for $\gamma\in\Gamma$ and $f\in C_c(G)$ with $\supp(f)\subseteq c^{-1}(\gamma)$ (a reduced coaction of a discrete group $\Gamma$ on a $C^*$-algebra $A$ is a $*$-homomorphism $\delta:A\to A\otimes C^*_r(\Gamma)$ (where $\otimes$ denotes the minimal tensor product, and $C^*_r(\Gamma)$ the reduced group $C^*$-algebra of $\Gamma$) that is nondegenerate in the sense that $\operatorname{span}\{\delta(a_1)(a_2\otimes f):a_1,a_2\in A,\ f\in C^*_r(\Gamma)\}$ is dense in $A\otimes C^*_r(\Gamma)$, and satisfies the identity $(\delta\otimes 1)\circ\delta=(1\otimes \delta_\Gamma)\circ\delta$ where $\delta_\Gamma:C^*_r(\Gamma)\to C^*_r(\Gamma)\otimes C^*_r(\Gamma)$ is the $*$-homomorphism that maps $\lambda_\gamma\to\lambda_\gamma\otimes\lambda_\gamma$ for each $\gamma\in\Gamma$, where $\gamma\mapsto\lambda_\gamma$ is the unitary representation of $\Gamma$ in $C^*_r(\Gamma)$).

If $c:G\to\Gamma$ is a cocycle, then $c^{-1}(e)$ (where $e$ is the identity of $\Gamma$) is an open subgroupoid of $G$. It follows that $C^*_r(c^{-1}(e))$ is a $C^*$-subalgebra of $C^*_r(G)$. We then have that $C^*_r(c^{-1}(e))$ is equal to the generalised fixed-point algebra $C^*_r(G)^{\delta_c}:=\{f\in C^*_r(G): \delta_c(f)=f\otimes e\}$ for $\delta_c$.

\chapter{Reconstruction of étale groupoids}\label{chapter:weyl}

We saw in the previous chapter that if $G$ is a locally compact Hausdorff étale groupoid, then we get a $C^*$-algebra $C^*_r(G)$ and an abelian $C^*$-algebra $C_0(\go)$ which contains an approximate identity of $C^*_r(G)$. 

We shall in this chapter see how it in some cases is possible to recover $G$ from $C^*_r(G)$ and $C_0(\go)$. This means that if $G$ and $H$ are locally compact Hausdorff étale groupoids satisfying certain conditions which we specify later, then there is a $*$-isomorphism from $C^*_r(G)$ to $C^*_r(H)$ that maps $C_0(\go)$ onto $C^*_r(H^{(0)})$ if and only if $G$ and $H$ are isomorphic as topological groupoids. 

If $G$ and $H$ are second countable locally compact Hausdorff étale topologically principal groupoids (meaining that $\iso(G)^\circ=\go$ and $\iso(H)^\circ=H^{(0)}$), then this follows from results in \cite{Re2008}. It was proven in the paper \cite{BrCaWh2017} that it also holds if $G$ and $H$ are the groupoids of countable directed graphs, even if $G$ and $H$ are not topologically free. 

If $G=\Z_4$ and $H=\Z_2\otimes\Z_2$, then $(C^*_r(G),C_0(\go))$ and $(C^*_r(H),C_0(H^{(0)}))$ are isomorphic, but $G$ and $H$ are not, so we cannot in general recover any locally compact Hausdorff étale groupoid from $C^*_r(G)$ and $C_0(\go)$. We shall see that if we assume that $\iso(G)^\circ$ is abelian and torsion-free, then this obstacle disappears. 

It is for some application interesting to recover not only $G$, but also a cocycle from $G$ to a discrete group $\Gamma$ from $C^*_r(G)$ and $C_0(\go)$. We shall see that this is also possible if we in addition to $C^*_r(G)$ and $C_0(\go)$ also include a coaction of $C^*_r(G)$.

We shall assume that $G$ is second countable which means that $C^*_r(G)$ is separable. Most of the result presentated in this chapter are taken from the  paper \cite{CaRuSiTo2017}.

\section{The extended Weyl groupoid of a weakly Cartan pair}

The following notion of weakly Cartan pairs was introduced in \cite{CaRuSiTo2017}. It is a generalisation of the notion of Cartan pairs introduced in \cite{Re2008}.

\begin{definition}\label{def:semicartan}
Let $A$ be a separable $C^*$-algebra. A $C^*$-subalgebra $D \subseteq A$ is called \emph{weakly Cartan} if
\begin{enumerate}
\item\label{it:sc abelian} $D$ is abelian,
\item\label{it:sc ai} $D$ contains an approximate identity for $A$,
\item\label{it:sc unital fibres} for each $\phi \in \widehat{D}$, the quotient
    $D'/J_\phi$ of $D'$ by the ideal $J_\phi := \ker(\phi)D'$ is a unital
    $C^*$-algebra, and
\item\label{it:sc unital sections} for each $\phi \in \widehat{D}$, there exist $d
    \in D$ and an open neighbourhood $U$ of $\phi$ such that $d + J_\psi =
    1_{D'/J_\psi}$ for all $\psi \in U$.
\end{enumerate}
We call the pair $(A, D)$ a \emph{weakly Cartan pair} of $C^*$-algebras.
\end{definition}

If $D \subseteq A$ is a Cartan subalgebra in the sense of Renault \cite{Re2008}, then it satisfies (\ref{it:sc abelian})~and~(\ref{it:sc ai}) by definition. Since Cartan subalgebras are by definition maximal abelian, we have $D' = D \cong C_0(\widehat{D})$, so each $J_\phi = \ker(\phi)$ and each $D'/J_\phi \cong \phi(D) = \C$. So $D$ satisfies \eqref{it:sc unital fibres}~and~\eqref{it:sc unital sections} by the Gelfand--Naimark theorem. That is, every Cartan subalgebra is in particular a weakly Cartan subalgebra, justifying our terminology.

It follows from the discusion in Chapter~\ref{chapter:groupoids} that if $G$ is a locally compact second-countable Hausdorff étale groupoid, then $(C^*_r(G),C_0(G^{(0)}))$ is a weakly Cartan pair. Also, if $c:G\to\Gamma$ is a cocycle from $G$ to a discrete group $\Gamma$, then $(C^*_r(c^{-1}(e)),C_0(G^{(0)}))$ is a weakly Cartan pair.

In is shown in \cite{CaRuSiTo2017} how one from a separable $C^*$-algebra $A$, a coaction $\delta$ of a countable discrete group $\Gamma$ on $A$, and a weakly Cartan subalgebra $D$ of the generalised fixed-point algebra $A^\delta$ can construct a second countable locally compact Hausdorff étale groupoid $\Hh(A,D,\delta)$ and a cocycle $c_\delta:\Hh(A,D,\delta)\to\Gamma$. 
For this we need the notion of a \emph{normaliser}. 

Normalisers were introduced in \cite{Ku1986a} and further studied in \cite{Re2008}.
Let $A$ be a separable $C^*$-algebra, and $D$ an abelian $C^*$-subalge\-bra of $A$ containing an approximate identity for $A$. A \emph{normaliser} $n$ of $D$ in $A$ is an element $n\in A$ such that $nDn^*\cup n^*Dn\subseteq D$. We write $N(D)$ for the collection of all normalisers of $D$ in $A$, and $\widehat{D}$ for the set of characters of $D$. For $d \in D$, we write $\osupp(d) :=\{\phi \in \widehat{D} : \phi(d) \not= 0\}$ and $I(d):=\{d'\in D : \osupp(d')\subseteq\osupp(d)\}$. Then $\osupp(d)$ is an open subset of $D'$, and $I(d)$ is an ideal of $D$. If $n\in N(D)$, then $nn^*,n^*n\in D$ and there is a homeomorphism $\alpha_n : \osupp(n^*n) \to \osupp(nn^*)$ such that $\phi(n^*n) \alpha_n(\phi)(d) = \phi(n^*dn)$ for all $d \in D$. We have that $\alpha_n^{-1}=\alpha_{n^*}$. Moreover, if $m,n\in N(D)$, then $mn \in N(D)$, $\osupp((mn)^*(mn)) =
    \alpha_n^{-1}(\osupp(m^*m)\cap\osupp(nn^*))$, and on this domain, $\alpha_m \circ \alpha_n =\alpha_{mn}$.


We assume now that $A$ is a separable $C^*$-algebra, $\Gamma$ is a discrete group, $\delta$ a is coaction of $\Gamma$ on $A$, and $D$ is a weakly Cartan subalgebra of the generalised fixed-point algebra $A^\delta$. 
We write
\[
D'_{A^\delta} := \{a \in A^\delta : ad = da\text{ for all } d \in D\}
\]
for the relative commutant of $D$ in $A^\delta$, and
\[
\pi_\phi : D'_{A^\delta} \to D'_{A^\delta}/J_\phi
\]
for the canonical quotient map.
It is proven in \cite[Corollary 1.6]{Qu1996} that $A^\delta$ contains an approximate identity for $A$. It follows that any approximate identity for of the
generalised fixed-point algebra $A^\delta$ is also an approximate identity for $A$, so $D$ contains an approximate identity for $A$. We say that a normaliser $n$ of $D$ is a \emph{homogeneous normaliser} if $n \in A_g$ for some $g \in \Gamma$. We write write $N_\star(D)$ for the collection of homogeneous normalisers of $A$. For $g \in \Gamma$ we write $N_g(D)$ for $N_\star(D) \cap A_g$.

The groupoid $\Hh(A,D,\delta)$ consists of equivalence classe of pairs $(n,\phi)$ where $n\in N_\star(D)$ and $\phi\in\osupp(n^*n)$. We need the following lemma from the paper \cite{CaRuSiTo2017} to define the equivalence relation used in the definition of $\Hh(A,D,\delta)$.

\begin{lemma}\label{lem:unitary-normalizers}
Let $n, m \in N(D)$ and $\phi \in \osupp(n^*n)\cap\osupp(m^*m)$, and suppose that there is an open neighbourhood $U$ of $\phi$ such that $U \subseteq \osupp(n^*n) \cap \osupp(m^*m)$ and $\alpha_n \vert_U = \alpha_m \vert_U$. Fix $d \in D$ with $\osupp(d) \subseteq U$ and $\phi(d) = 1$, and let 
\begin{equation}\label{eq:wdef}
w := \phi(n^*n)^{-1/2}\phi(m^*m)^{-1/2}d n^*m d.
\end{equation}
Then $w\in D'_{A^\delta}$ and $\pi_\phi(w)$ is unitary in $D'_{A^\delta}/J_\phi$. Moreover, $\pi_\phi(w)$ is independent of the choices of $U$ and $d$.
\end{lemma}

We write
\[
U^\phi_{n^*m} := \pi_\phi(w)
\]
for any $w$ of the form~\eqref{eq:wdef}. If $\phi$ is clear from context, we just write $U_{n^*m}$ for $U^\phi_{n^*m}$. It is easy to check that if $n,m,n' \in N(D)$, $\phi\in\osupp(n^*n)\cap\osupp(m^*m)\cap\osupp((n')^*n')$, and there is an open neighbourhood $U$ of $\phi$ such that $U \subseteq$\linebreak $\osupp(n^*n) \cap \osupp(m^*m)$ and $\alpha_m|_U = \alpha_n|_U$, and an open neighbourhood $U'$ of $\phi$ such that $U' \subseteq \osupp( (n')^*n') \cap \osupp(m^*m)$ and $\alpha_m|_U = \alpha_{n'}|_U$, then
\begin{enumerate}
	\item\label{it:u1} $U_{n^*n}=1_{D'_{A^\delta}/J_\phi}$,
	\item\label{it:u2} $U_{n^*m}^*=U_{mn^*}$, and
	\item\label{it:u3} $U_{n^*m}U_{m^*n'}=U_{n^*n'}$.
\end{enumerate}

We are now ready to define the equivalence relation used in the definition of $\Hh(A,D,\delta)$.
Define a relation on $\{(n,\phi) :
n \in N_\star(D),\ \phi\in\osupp(n^*n)\}$ by $(n, \phi) \sim (m,\psi)$ if and only if
\begin{itemize}
\item[(R1)] $\phi = \psi$,
\item[(R2)] $n^*m \in A^\delta$,
\item[(R3)] there exists an open neighbourhood $U$ of $\phi$ in $\widehat{D}$ such
    that $U\subseteq \osupp(n^*n)\cap\osupp (m^*m)$ and $\alpha_m \vert_U = \alpha_n \vert_U$, and
\item[(R4)] the unitary $U^\phi_{n^*m}$ belongs to the
    connected component $\mathcal{U}_0(D'_{A^\delta}/J_\phi)$ of the identity in the
    unitary group of $D'_{A^\delta}/J_\phi$.
\end{itemize}
Then $\sim$ is an equivalence relation.

The following 4 results can be found in \cite{CaRuSiTo2017}.

\begin{theorem}\label{thm:lche gpd}
Let $A$ be a separable $C^*$-algebra, $\delta$ a coaction of a countable discrete group
$\Gamma$ on $A$, and $D$ a weakly Cartan subalgebra of $A^\delta$. For $n \in N_\star(D)$ and
$\phi \in \supp(n^*n)$, let $[n,\phi]$ denote the equivalence class
of $(n,\phi)$ under $\sim$. Define
\[
\Hh(A,D,\delta) := \{[n, \phi] : n \in N_\star(D),\ \phi \in \supp(n^*n)\}.
\]
There are maps
\begin{gather*}
    r, s \colon \Hh(A, D, \delta) \to \widehat{D},\\
       M \colon \Hh(A, D, \delta) \mathbin{_s\times_r} \Hh(A, D, \delta) \to \Hh(A, D, \delta),\text{ and}\\
       I \colon \Hh(A, D, \delta) \to \Hh(A, D, \delta)
\end{gather*}
such that
\begin{gather*}
    r([n, \phi] ) = \alpha_n(\phi), \quad s([n,\phi]) = \phi, \\
    \quad M([n, \phi],[m,\psi]) = [nm, \psi], \quad \text{and} \quad I([n,\phi]) = [n^*, \alpha_n(\phi)].
\end{gather*}
Moreover, $\Hh(A, D, \delta)$ is a groupoid under these operations, and there is a cocycle $c_\delta : \Hh(A, D, \delta) \to \Gamma$ such that $c_\delta([n,\phi]) = g$ if and only if $n \in A_g$.

For $n \in N_\star(D)$ and an open set $X \subseteq \widehat{D}$ contained in $\supp(n^*n)$ let
\[
    Z(n, X) := \{[n,\phi] : \phi \in X\} \subseteq \Hh(A, D, \delta).
\]
Then 
\begin{equation}\label{eq:basis}
    \{Z(n, X) : n \in N_\star(D),\ X \subseteq \widehat{D}\text{ is open and } X\subseteq\supp(n^*n)\}
\end{equation}
constitutes a basis for a second countable locally compact Hausdorff \'etale topology on the groupoid $\Hh(A, D, \delta)$. The cocycle $c_\delta$ is continuous with respect to this topology.
\end{theorem}

Let $\Gamma$ be a discrete group with identity $e$. If $G$ a locally compact Hausdorff \'etale groupoid and $c : G \to \Gamma$ is a cocycle, then $c^{-1}(e)$ is open subgroupoid of $G$ and $C^*_r(c^{-1}(e))$ is equal to the generalised fixed-point algebra of the coaction $\delta_c$ of $\Gamma$ on $C_r^*(G)$. It follows that $C_0(\go)$ is a weakly Cartan subalgebra of $C^*_r(c^{-1}(e))$.
If $x\in\go$, then we denote by $\widehat{x}$ the character of $C_0(\go)$ given by $\widehat{f}=f(x)$.  

\begin{theorem}\label{prp:isomorphism}
Let $\Gamma$ be a discrete group, $G$ a locally compact Hausdorff \'etale groupoid and $c : G \to \Gamma$ a cocycle. Suppose that $\Io{c^{-1}(e)}$ is torsion-free and abelian.
Then there is an isomorphism $\theta : G \to \Hh(C^*_r(G), C_0(\go), \delta_c)$ such that $c_{\delta_c}\circ\theta=c$ and 
such that for $\gamma \in G$ and any $n \in C_c(G)$ supported on a
bisection contained in $c^{-1}(c(\gamma))$ and satisfying $n(\gamma) = 1$, we have $\theta(\gamma) = [n,
\widehat{s(\gamma)}]$.
\end{theorem}

By letting $\Gamma$ be the trivial group and $e$ the trivial cocycle on $G$, we obtain the following corollary.

\begin{corollary}
Let $G$ a locally compact Hausdorff \'etale groupoid such that $\Io{G}$ is torsion-free and abelian.
Then there is an isomorphism $\theta : G \to \Hh(C^*_r(G), C_0(\go), \delta_e)$
such that for $\gamma \in G$ and any $n \in C_c(G)$ supported on a
bisection and satisfying $n(\gamma) = 1$, we have $\theta(\gamma) = [n,
\widehat{s(\gamma)}]$.
\end{corollary}

We also obtain from Theorem~\ref{prp:isomorphism} the following $C^*$-rigidity result.

\begin{theorem}\label{thm:1}
Let $\Gamma$ be a discrete group, and let $(G_1,c_1), (G_2,c_2)$ be $\Gamma$-graded
locally compact Hausdorff \'etale groupoids such that each $\Io{c_i^{-1}(e)}$ is
torsion-free and abelian.
\begin{enumerate}
\item Suppose that $\kappa : G_2 \to G_1$ is an isomorphism satisfying
    $c_1\circ\kappa=c_2$. Then there is an isomorphism $\phi : C^*_r(G_1)\to
    C^*_r(G_2)$ such that $\phi(f) = f \circ \kappa$ for $f \in C_c(G_1)$. We have
    $\phi(C_0(G_1^{(0)})) = C_0(G_2^{(0)})$ and $\delta_{c_2}\circ\phi =
    (\phi\otimes\id)\circ\delta_{c_1}$.
\item Suppose that $\phi:C^*_r(G_1)\to C^*_r(G_2)$ is an isomorphism satisfying\linebreak
    $\phi(C_0(G_1^{(0)}))=C_0(G_2^{(0)})$ and
    $\delta_{c_2}\circ\phi=(\phi\otimes\id)\circ\delta_{c_1}$. Then there is an
    isomorphism $\kappa : G_2 \to G_1$ such that $\kappa|_{G_2^{(0)}}$ is the
    homeomorphism induced by $\phi|_{C_0(G_1^{(0)})}$ and $c_1 \circ \kappa = c_2$.
\end{enumerate}
\end{theorem}


\chapter{$C^*$-rigidity of dynamical systems}\label{chapter:dyn}
We shall now look at how to obtain $C^*$-rigidity results of some dynamical systems from Theorem~\ref{thm:1}. 

\section{Homeomorphisms of compact Hausdorff spaces}

We begin by looking at homeomorphisms of compact Hausdorff spaces. The following theorem is proven in \cite{CaRuSiTo2017}.

\begin{theorem}
Suppose that $\sigma : X \to X$ and $\tau : Y \to Y$ are homeomorphisms of compact
Hausdorff spaces. The following are equivalent:
\begin{enumerate}
    \item $G(X,\sigma)$ and $G(Y,\tau)$ are isomorphic;
    \item $C(X)\times_\sigma\Z \cong C(Y)\times_\tau\Z$ via an isomorphism that
        maps $C(X)$ to $C(Y)$; and
    \item there exist decompositions $X=X_1\sqcup X_2$ and $Y=Y_1\sqcup Y_2$ into
        disjoint open invariant sets such that $\sigma|_{X_1}$ is conjugate to $\tau|_{Y_1}$ and
        $\sigma|_{X_2}$ is conjugate to $\tau^{-1}|_{Y_2}$.
\end{enumerate}
If $\sigma$ and $\tau$ are topologically transitive or $X$ and $Y$ are connected, then
these conditions hold if and only $\sigma$ and $\tau$ are flip-conjugate.
\end{theorem}

Since $C^*_r(G(X,\sigma))$ is isomorphic to $C(X)\times_\sigma\Z$ and $C^*_r(G(Y,\tau))$ is isomorphic to $C(Y)\times_\tau\Z$, the equivalence of 1. and 2. follows from Theorem~\ref{thm:1}. It is straightforward to check that 3. implies 1. For the proof of $1.\implies 3.$ one uses that if $\theta:G(X,\sigma)\to G(Y,\tau)$ is and isomorphism, then we get a homeomorphism $h : X \to Y$ by letting
$\theta(x,0,x)=(h(x),0,h(x))$ for $x\in X$, and a map $f : X \times \Z \to \Z$ by letting
$f(n,x):=c_Y(\theta(x,n,\sigma^n(x)))$ where $c_X : G(X, \sigma) \to \Z$ and $c_Y : G(Y, \tau) \to
\Z$ are the canonical cocycles, such that
\begin{equation*}
f(m+n,x) = f(m,x) + f(n,\sigma^m(x))
\end{equation*}
for all $m,n\in\Z$ and $x\in X$,
and then $f(\cdot,x)$ is for $x\in X$ a bijection of $\Z$ with inverse $n\mapsto
c_X(\theta^{-1}(h(x),n,\tau^n(h(x))))$. There is then a positiv integer $N$ such that 
\begin{align*}
X_1 &:= \{x\in X : f(n,x)>0\text{ and }f(-n,x) < 0\text{ for }n> N\}\quad\text{ and}\\
X_2 &:= \{x\in X : f(n,x)<0\text{ and }f(-n,x) > 0\text{ for }n> N\}
\end{align*}
are clopen $\sigma$-invariant subsets such that $X=X_1\sqcup X_2$ and $\sigma|_{X_1}$ is conjugate to $\tau|_{h(X_1)}$ and $\sigma|_{X_2}$ is conjugate to $\tau^{-1}|_{h(X_2)}$.

\section{Group actions}

Next, we look at group actions.

Let $\Gamma$ be a discrete group acting on the right of a second-countable locally compact Hausdorff space $X$. 

For $x\in X$, we write
\[
\Stab(x) := \{\gamma\in\Gamma : x\gamma=x\}
\]
for the \emph{stabilizer subgroup} of $x$ in $\Gamma$; observe that then $(X \rtimes \Gamma)^x_x = \{x\} \times \Stab(x)$. We also consider the \emph{essential stabilizer subgroup}
\[
\Stab^{\ess}(x) := \{\gamma \in \Gamma : \gamma \in \Stab(y)\text{ for all $y$ in some neighbourhood $U$ of $x$}\}.
\]
A Baire-category argument shows that $(X, \Gamma)$ is topologically free if and only if $\overline{\{x \in X : \Stab(x) = \{e\}\}} = X$.

\begin{definition}\label{dfn:coe groups}
Let $\Gamma$ and $\Lambda$ be countable discrete groups acting on second-countable locally compact Hausdorff spaces $X$ and $Y$. A \emph{continuous orbit equivalence} between $(X,\Gamma)$ and $(Y,\Lambda)$ is a triple $(h,\phi,\eta)$ consisting of a homeomorphism $h : X \to Y$, and continuous functions $\phi : X\times \Gamma \to \Lambda$ and $\eta : Y\times\Lambda \to \Gamma$ such that $h(x\gamma) = h(x)\phi(x,\gamma)$ for all $x,\gamma$ and $h^{-1}(y\lambda)=h^{-1}(y)\eta(y,\lambda)$ for all $y, \lambda$. We call $h$ the \emph{underlying homeomorphism} of $(h,\phi,\eta)$.
\end{definition}

For topologically free systems, the intertwining condition appearing in Definition~\ref{dfn:coe groups} has some important consequences.

\begin{lemma}\label{lem:coe consequences}
Let $\Gamma \curvearrowright X$ and $\Lambda \curvearrowright$ be topologically free actions of countable discrete groups on second-countable locally compact Hausdorff spaces. Let $(h, \phi, \eta)$ be a continuous orbit equivalence from $(X, \Gamma)$ to $(Y, \Lambda)$. Fix $x \in X$. We have
\[
    \phi(x, \gamma\gamma') = \phi(x,\gamma)\phi(x\gamma, \gamma')\quad\text{for all $\gamma,\gamma' \in \Gamma$,}
\]
the maps $\theta_x := \phi(x, \cdot) : \Gamma \to \Lambda$ is a bijection that carries $e_\Gamma$ to $e_\Lambda$ and restricts to bijections $\Stab(x) \to \Stab(h(x))$ and $\Stab^{\ess}(x) \to \Stab^{\ess}(h(x))$.
\end{lemma}

For general systems, the conclusions of Lemma~\ref{lem:coe consequences} are not automatic features of continuous orbit equivalences, prompting the following definition.

\begin{definition}
Let $\Gamma \curvearrowright X$ and $\Lambda \curvearrowright Y$ be actions of countable discrete groups on second-countable locally compact Hausdorff spaces. Consider a map $\phi : X \times \Gamma \to \Lambda$.
\begin{enumerate}
\item We call $\phi$ a \emph{cocycle} if $\phi(x, \gamma\gamma') =
    \phi(x,\gamma)\phi(x\gamma, \gamma')$ for all $x, \gamma, \gamma'$.
\item Let $h : X \to Y$ be a homeomorphism. We say that $(h,\phi)$ \emph{preserves
    stabilizers} if $\phi(x,\cdot)$ restricts to bijections $\Stab(x) \to
    \Stab(h(x))$, and that $(h,\phi)$ \emph{preserves essential stabilizers} if
    $\phi(\cdot, x)$ restricts to bijections $\Stab^{\ess}(x) \to
    \Stab^{\ess}(h(x))$.
\end{enumerate}
\end{definition}

The following proposition then follows from Theorem~\ref{thm:1}.

\begin{proposition} \label{prp:group action}
Let $\Gamma \curvearrowright X$ and $\Lambda \curvearrowright Y$ be actions of countable discrete groups on second-countable locally compact Hausdorff spaces. Suppose that $h : X \to Y$ is a homeomorphism and $\phi : X\times \Gamma \to \Lambda$ is continuous. The following are equivalent:
\begin{enumerate}
    \item there is an isomorphism $\Theta : X \rtimes \Gamma \to Y \rtimes \Lambda$ such that $\Theta(x,e)=(h(x),e)$ and $\Theta(x,\gamma) = (h(x),\phi(x,\gamma))$ for all $x\in X$ and $\gamma\in\Gamma$;
    \item $\phi$ is a cocycle, $(h,\phi)$ preserves stabilizers, and there is a map $\eta : Y \times \Lambda \to \Gamma$ such that $(h,\phi,\eta)$ is a continuous orbit equivalence; and
    \item $\phi$ is a cocycle, $(h,\phi)$ preserves essential stabilizers, and there is a map $\eta : Y \times \Lambda \to \Gamma$ such that $(h,\phi,\eta)$ is a continuous orbit equivalence.
\end{enumerate}
\end{proposition}

This gives the following generalisation of Li's rigidity theorem~\ref{thm:group}.

\begin{corollary}
Let $\Gamma \curvearrowright X$ and $\Lambda \curvearrowright Y$ be actions of countable discrete groups on second-countable locally compact Hausdorff spaces. 

Suppose that $\Stab^{\ess}(x)$ and $\Stab^{\ess}(y)$ are torsion-free and abelian for all $x\in X$ and all $y\in Y$, and that $h$ is a homeomorphism from $X$ to $Y$. The following are equivalent.
\begin{enumerate}
    \item There exist cocycles $\phi : X \times \Gamma \to \Lambda$ and $\eta : Y \times \Lambda \to \Gamma$ such that $(h, \phi, \eta)$ is a continuous orbit equivalence from $(X,\Gamma)$ to $(Y,\Lambda)$ and $(h, \phi)$ and $(h^{-1}, \eta)$ preserve essential stabilizers.
    \item There is an isomorphism $\Theta : X \rtimes \Gamma \to Y \rtimes \Lambda$ such that $\Theta(x,e)=(h(x),e)$ for $x\in X$.
    \item There is a $*$-isomorphism $\phi:C_0(X)\rtimes_r\Gamma\to C_0(Y)\rtimes_r\Lambda$\linebreak such that $\phi(C_0(X))=C_0(Y)$ and $\phi(f)=f\circ h^{-1}$ for $f\in C_0(X)$.
\end{enumerate}
\end{corollary}

\section{Local homeomorphisms}

Let us finally look at how Theorem~\ref{thm:1} applied to $G(X,\sigma)$ where $X$ is a locally compact Hausdorff space, and $\sigma$ is a local homeomorphism from an open subset of $X$ to an open subset of $X$, can be used to prove Theorem~\ref{thm:ck}.

Suppose $X$ is a compact Hausdorff space and $\sigma:X\to X$ is a surjective local homeomorphism. Then there is a continuous map $T_{(X,\sigma)}:G(X,\sigma)\to G(X,\sigma)$ such that $T_{(X,\sigma)}(x,k,y)=(\sigma(x),k,\sigma(y))$. One can show that there is a continuous positive map $\phi_{(X,\sigma)}:C^*_r(G(X,\sigma))\to C^*_r(G(X,\sigma))$ such that $\phi_{(X,\sigma)}(f)(\eta)=f(T_{(X,\sigma)}(\eta))$ for $f\in C_c(G(X,\sigma))$ and $\eta\in G(X,\tau)$, and that the following theorem holds.
\begin{theorem}\label{thm:pos}
Let $X$ and $Y$ be compact Hausdorff spaces and $\sigma:X\to X$ and $\tau:Y\to Y$ surjective local homeomorphisms. The following are equivalent.
\begin{enumerate}
	\item $(X,\sigma)$ and $(Y,\tau)$ are conjugate.
	\item There is an isomorphism $\theta:G(X,\sigma)\to G(Y,\tau)$ such that $\theta\circ T_{(X,\sigma)}=T_{(Y,\tau)}\circ\theta$.
	\item There is a $*$-isomorphism $\psi:C^*_r(G(X,\sigma))\to C^*_r(G(Y,\tau))$ such that $\psi(C(X))=C(Y)$ and $\psi\circ \phi_{(X,\sigma)}=\phi_{(Y,\tau)}\circ\psi$.
\end{enumerate}
\end{theorem}
 
Statement 1. in Theorem~\ref{thm:ck} then follows from Theorem~\ref{thm:pos}.

Let $X$ be a locally compact Hausdorff space, and let $\sigma : U \to V$ be a local homeomorphism from an open set $U \subseteq X$ to an open set $V \subseteq X$. For $x\in X$ we define the \emph{stabiliser group} at $x$ by
\begin{equation*}
\Stab(x) := \{m-n:m,n\in\N_0,\ x\in U_m\cap U_n,\text{ and }\sigma^n(x)=\sigma^m(x)\} \subseteq \Z,
\end{equation*}
and we define the \emph{essential stabiliser group} at $x$ by
\begin{multline*}
\Stab^{\ess}(x) := \{m-n:m,n\in\N_0\text{ and there is an open neighbourhood }\\
    O\subseteq U_m\cap U_n\text{ of }x\text{ such that }\sigma^n|_O = \sigma^m|_O\} \subseteq \Stab(x).
\end{multline*}
With the convention that $\min(\emptyset) = \infty$, we define the \emph{minimal
stabiliser} of $x$ to be
\[
\Stab_{\min}(x) := \min\{n \in \Stab(x) : n \ge 1\},
\]
and the \emph{minimal essential stabiliser} at $x$ to be
\[
    \Stab^{\ess}_{\min}(x) := \{n \in \Stab^{\ess}(x) : n \ge 1\}.
\]

The isotropy subgroupoid of
$\G(X,\sigma)$ is $\{(x,n,x):x\in X,\ n\in\Stab(x)\}$, and the interior of the isotropy is $\Io{\G(X,\sigma)}=\{(x,n,x):x\in X,\ n\in\Stab^{\ess}(x)\}$. So $\Io{\G(X,\sigma)}$ is torsion-free and abelian. Taking $\Gamma = \{e\}$ and $c : \G \to \Gamma$ the trivial cocycle, we obtain a (trivially) graded groupoid $\G$.

\begin{definition}
Let $X$ and $Y$ be locally compact Hausdorff spaces, $\sigma$ a local homeomorphism from an open subset $U_\sigma$ of $X$ to an open subset $V_\sigma$ of $X$, and $\tau$ a local homeomorphism from an open subset $U_\tau$ of $Y$ to an open subset $V_\tau$ of $Y$. We say that $(X, \sigma)$ and $(Y, \tau)$ are \emph{continuous orbit equivalent} if there exist a homeomorphism $h : X \to Y$  and continuous maps $k,l : U_\sigma \to \N_0$ and $k', l' : U_\tau \to \N_0$ such that
\[
\tau^{l(x)}(h(x)) = \tau^{k(x)}(h(\sigma(x)))
    \quad\text{ and }\quad
\sigma^{l'(y)}(h^{-1}(y)) =
\sigma^{k'(y)}(h^{-1}(\tau(y)))
\]
for all $x, y$. We call $(h,l,k,l',k')$ a continuous orbit equivalence and we call $h$ the \emph{underlying homeomorphism}. We say that $(h,l,k,l',k')$ \emph{preserves stabilisers} if $\Stab_{\min}(h(x))  <\infty \iff \Stab_{\min}(x) <\infty$, and
\begin{align*}
&\bigg|\sum_{n=0}^{\Stab_{\min}(x)-1} l(\sigma^n(x)) - k(\sigma^n(x)) \bigg| =
\Stab_{\min}(h(x))
\text{ and}\\
&\bigg|\sum_{n=0}^{\Stab_{\min}(y)-1} l'(\tau^n(y)) - k'(\tau^n(y)) \bigg| =
\Stab_{\min}(h^{-1}(y))
\end{align*}
whenever $\Stab(x), \Stab(y)$ are nontrivial, $\sigma^{\Stab_{\min}(x)}(x)=x$, and
$\tau^{\Stab_{\min}(y)}(y)=y$.

Likewise, we say that $(h,l,k,l',k')$ \emph{preserves essential stabilisers} if\linebreak
$\Stab^{\ess}_{\min}(h(x)) < \infty\iff \Stab^{\ess}_{\min}(x) < \infty$, and
\begin{align*}
&\bigg|\sum_{n=0}^{\Stab^{\ess}_{\min}(x)-1} \left( l(\sigma^n(x)) - k(\sigma^n(x)) \right) \bigg| =
\Stab^{\ess}_{\min}(h(x))
\text{ and }\\
&\bigg|\sum_{n=0}^{\Stab^{\ess}_{\min}(y)-1}\left( l'(\tau^n(y)) - k'(\tau^n(y)) \right) \bigg| =
\Stab^{\ess}_{\min}(h^{-1}(y))
\end{align*}
whenever $\Stab^{\ess}_{\min}(x), \Stab^{\ess}_{\min}(y) < \infty$,
$\sigma^{\Stab^{\ess}_{\min}(x)}(x)=x$, and $\tau^{\Stab^{\ess}_{\min}(y)}(y)=y$.
\end{definition}

By using Theorem~\ref{thm:1} the following theorem is shown in \cite{CaRuSiTo2017}.

\begin{theorem}\label{thm:DR oe<->gi}
Let $X$ and $Y$ be locally compact Hausdorff spaces, $\sigma$ a local homeomorphism from an open subset $U_\sigma$ of $X$ to an open subset $V_\sigma$ of $X$, and $\tau$ a local homeomorphism from an open subset $U_\tau$ of $Y$ to an open subset $V_\tau$ of $Y$, and suppose that $h : X
\to Y$ is a homeomorphism. Then the following are equivalent:
\begin{enumerate}
\item\label{it:DR oe} there is a stabiliser-preserving continuous orbit equivalence
    from $(X, \sigma)$ to $(Y, \tau)$ with underlying homeomorphism $h$;
    \item there is a essential stabiliser-preserving continuous orbit equivalence
    from $(X, \sigma)$ to $(Y, \tau)$ with underlying homeomorphism $h$;
\item\label{it:DR G iso} there is a groupoid isomorphism $\Theta : \G(X, \sigma) \to
    \G(Y, \tau)$ such that $\Theta|_X = h$; and
\item\label{it:DR Cstar iso} there is an isomorphism $\phi : C^*(\G(X, \sigma)) \to
    C^*(\G(Y, \tau))$ such that\linebreak $\phi(C_0(X)) = C_0(Y)$ with $\phi(f) = f\circ
    h^{-1}$ for $f \in C_0(Y)$.
\end{enumerate}
\end{theorem}

Statement 3. in Theorem~\ref{thm:ck} then follows from Theorem~\ref{thm:DR oe<->gi} and\linebreak \cite[Corollary 4.6]{CaWi2016}.

\begin{definition}
Let $X$ and $Y$ be locally compact Hausdorff spaces, $\sigma$ a local homeomorphism from an open subset $U_\sigma$ of $X$ to an open subset $V_\sigma$ of $X$, and $\tau$ a local homeomorphism from an open subset $U_\tau$ of $Y$ to an open subset $V_\tau$ of $Y$. We say that $(X, \sigma)$
and $(Y, \tau)$ are \emph{eventually conjugate} if there is a stabiliser-preserving
continuous orbit equivalence $(h,l,k,l',k')$ from $(X, \sigma)$ to $(Y, \tau)$ such that
$l(x) = k(x) + 1$ for all $x \in X$.
\end{definition}

Given $(X, \sigma)$, there is an action $\gamma^X : \T \to \Aut(C^*(\G(X, \sigma)))$
such that\linebreak $\gamma^X_z(f)(x, n, x') = z^n f(x, n, x')$ for all $z \in \T$, $(x, n, x')
\in \G(X, \sigma)$ and $f \in C_c(\G(X, \sigma))$.

By applying Theorem~\ref{thm:1} to $G(X,\sigma)$ and the coaction $c_X:G(X,\sigma)\to\Z$ given by $c_X(x,k,y)=k$, the following theorem is shown in \cite{CaRuSiTo2017}.

\begin{theorem}\label{thm:ev conj equivalences}
Let $X$ and $Y$ be locally compact Hausdorff spaces, $\sigma$ a local homeomorphism from an open subset $U_\sigma$ of $X$ to an open subset $V_\sigma$ of $X$, and $\tau$ a local homeomorphism from an open subset $U_\tau$ of $Y$ to an open subset $V_\tau$ of $Y$. Then the following are equivalent:
\begin{enumerate}
\item\label{it:DR evconj} there is an eventual conjugacy from $(X, \sigma)$ to $(Y,
    \tau)$ with underlying homeomorphism $h$;
\item\label{it:DR G c-iso} there is an isomorphism $\Theta : \G(X, \sigma) \to \G(Y,
    \tau)$ such that $\Theta|_X = h$ and $c_X = c_Y \circ \Theta$; and
\item\label{it:DR graded iso} there is an isomorphism $\phi : C^*(\G(X, \sigma)) \to
    C^*(\G(Y, \tau))$ such that\linebreak $\phi(C_0(X)) = C_0(Y)$, with $\phi(f) = f\circ
    h^{-1}$ for $f \in C_0(X)$, and $\phi \circ \gamma^X_z = \gamma^Y_z \circ \phi$.
\end{enumerate}
\end{theorem}

Statement 2. in Theorem~\ref{thm:ck} then follows from Theorem~\ref{thm:ev conj equivalences}.

Groupoids $\G_1$ and $\G_2$ are \emph{equivalent} if there is a topological space $Z$
carrying commuting free and proper actions of $\G_1$ and $\G_2$ on the left and right
respectively such that $r : Z \to \go_1$ and $s : Z \to \go_2$ induce homeomorphisms
$\G_1\backslash Z \cong \go_2$ and $Z/\G_2 \cong \go_1$. If $\Gamma$ is a discrete group and $c_1:G_1\to \Gamma$ and $c_2:G_2\to \Gamma$ are cocycles, then a \emph{graded} $(\G_1, c_1)$--$(\G_2, c_2)$-equivalence consists of a
$\G_1$--$\G_2$-equivalence $Z$ and a continuous map $c_Z : Z \to \Gamma$ satisfying
$c_Z(\gamma \cdot z \cdot \eta) = c_1(\gamma)c_Z(z)c_2(\eta)$ for all $\gamma, z, \eta$.

Suppose that for $i = 1,2$, $A_i$ is a $C^*$-algebra carrying a coaction $\delta_i$ of a
discrete group $\Gamma$, and $D_i \subseteq A_i^{\delta_i}$ is a subalgebra. We say that
$(A_1, D_1)$ and $(A_2, D_2)$ are \emph{equivariantly Morita equivalent} if there are an
$A_1$--$A_2$-imprimitivity bimodule $X$ and a right-Hilbert bimodule morphism $\zeta: X
\to M(X\otimes C^*_r(\Gamma))$ such that $(\zeta\otimes\id_\Gamma)\circ\zeta =
(\id_X\otimes\delta_\Gamma)\circ\zeta$ and for each $g \in \Gamma$, the subspace $X_g :=
\{x \in X : \zeta(x) = x \otimes \lambda_g\}$ satisfies
\begin{equation*}
X_g = \clsp\{\xi \in X_g : {_{A_1}\langle \xi, \xi\cdot D_2\rangle} \subseteq D_1
    \text{ and }\langle D_1 \cdot \xi, \xi\rangle_{A_2} \subseteq D_2\}.
\end{equation*}
If $\Gamma = \{0\}$, we say that $(A_1, D_1)$ and $(A_2, D_2)$ are \emph{Morita
equivalent}.

We let $\R$ denote the discrete groupoid obtained from the equivalence relation $\sim$ on $\NN$ defined by $m\sim n$ for all $m,n\in\NN$.

Let $X$ be a locally compact Hausdorff space, and $\sigma : X \to X$ a local homeomorphism. Let $\widetilde{X} := X \times \NN$ with the product topology and define a (surjective) local homeomorphism $\tilde{\sigma} : \widetilde{X} \to \widetilde{X}$ by $\tilde{\sigma}(x,0) = (\sigma(x), 0)$ and $\tilde{\sigma}(x, n+1) = (x,n)$. We call $(\tilde{X}, \tilde{\sigma})$ the \emph{stabilisation} of $(X, \sigma)$. Then $\G(\widetilde{X}, \tilde\sigma) \cong \G(X, \sigma) \times \R$ via $((x,m), p, (y,n)) \mapsto \big((x, p-m+n,y), (m,n)\big)$.

The following theorem is shown in \cite{CaRuSiTo2017}.

\begin{theorem}\label{thm:mor}
Let $\sigma : X \to X$ and $\tau : Y \to Y$, be local homeomorphisms of second-countable
locally compact totally disconnected Hausdorff spaces. The following are equivalent:
    \begin{enumerate}
    \item there is a stabiliser-preserving continuous orbit equivalence from
        $(\widetilde{X}, \tilde\sigma)$ to $(\widetilde{Y}, \tilde{\tau})$;
    \item $G(X, \sigma)$ and $G(Y, \tau)$ are equivalent;
    \item $(C^*(\G(X, \sigma)), C_0(X))$ and $(C^*(\G(Y, \tau)), C_0(Y))$ are Morita
        equivalent; and
    \item there is an isomorphism $C^*(\G(X, \sigma)) \otimes \Kk \to C^*(\G(Y,
        \tau)) \otimes \Kk$ that carries $C_0(X) \otimes \Cc$ to $C_0(Y) \otimes
        \Cc$.
    \end{enumerate}
\end{theorem}

Statement 5. in Theorem~\ref{thm:ck} can be shown by using Theorem~\ref{thm:mor} and results from \cite{CaRuSi2017}.

Let $\sigma:X\to X$ be a surjective local homeomorphism of a compact Hausdorff space. Let $\overline{X} := \{\xi \in X^{\Z}:\sigma(\xi_n)=\xi_{n+1}\text{ for every }n\in\Z\}$, and define $\overline{\sigma}:\overline{X}\to\overline{X}$ by $\overline{\sigma}(\xi)_n = \sigma(\xi_n)$. We call $\sigma$ \emph{expansive} if there is a metrisation $(X,d)$ of $X$ and an $\epsilon>0$ such that $\sup_n d(\sigma^n(x),\sigma^n(x')) < \epsilon \implies x = x'$. We call $\epsilon$ an \emph{expansive constant for $(X, d, \sigma)$}.

The following theorem is shown in \cite{CaRuSiTo2017}.

\begin{theorem}\label{cor:two-sided conjugacy}
Let $X, Y$ be second-countable locally compact totally disconnected Hausdorff spaces and let $\sigma : U_\sigma \to V_\sigma$ and $\tau : U_\tau \to V_\tau$ be local homeomorphisms between open subsets of $X,Y$. The following are equivalent:
    \begin{enumerate}
    \item there are continuous, open maps $f:X\to Y$ and $f':Y\to X$, and continuous
        maps $a:X\to\NN$, $k:U_\sigma\to\NN$, $a':Y\to\NN$, and
        $k':U_\tau\to\NN$ such that $\sigma^{a(x)}(f'(f(x)))=\sigma^{a(x)}(x)$
        for $x\in X$, $\tau^{k(x)}(f(\sigma(x)))=\tau^{k(x)+1}(f(x))$ for $x\in
        U)\sigma$, $\tau^{a'(y)}(f(f'(y)))=\tau^{a'(y)}(y)$ for $y\in Y$, and\linebreak
        $\sigma^{k'(y)}(f'(\tau(y)))=\sigma^{k'(y)+1}(f'(y))$ for $y\in U_\tau$;
    \item there is a graded $(G(X,\sigma),c_X)$--$(G(Y,\tau),c_Y)$-equivalence $(Z,
        c_Z)$ such that $c_Z^{-1}(0)$ is a $c_X^{-1}(0)$--$c_Y^{-1}(0)$-equivalence;
    \item there is an isomorphism $\kappa :G(X,\sigma)\times\R\to G(Y,\tau)\times\R$
        such that $\overline{c}_Y\circ\kappa=\overline{c}_X$;
    \item there is a $\T$-equivariant Morita equivalence between $(C^*(\G(X,
        \sigma)), C_0(X))$ and $(C^*(\G(Y, \tau)), C_0(Y))$ with respect to the gauge
        actions $\gamma^X$ and $\gamma^Y$ whose fixed-point submodule is a $C^*(\G(X,
        \sigma)$--$C^*(\G(Y, \tau)$-imprimitivity bimodule; and
    \item there is an isomorphism $C^*(\G(X, \sigma)) \otimes \Kk \to C^*(\G(Y,
        \tau)) \otimes \Kk$ that carries $C_0(X) \otimes \Cc$ to $C_0(Y) \otimes
        \Cc$, and intertwines the actions $\gamma^X \otimes \id$ and $\gamma^Y
        \otimes \id$.
    \end{enumerate}
    If $X$ and $Y$ are compact, $U_\sigma=V_\sigma=X$, and
    $U_\tau=V_\tau=Y$, then each of the above five conditions implies
    \begin{itemize}
        \item[6.] $(\overline{X},\overline{\sigma})$ and $(\overline{Y},\overline{\tau})$ are conjugate.
    \end{itemize}
    If in addition $\sigma$ and $\tau$ are expansive, then all six conditions are
    equivalent.
\end{theorem}

Statement 4. in Theorem~\ref{thm:ck} then follows from Theorem~\ref{cor:two-sided conjugacy}.


\begin{thebibliography}{}

\bibitem[Arklint et~al., 2017]{AREIRU2017}
Arklint, S.~E., Eilers, S., \& Ruiz, E. (2017).
\newblock A dynamical characterization of diagonal-preserving ${\ast }$-isomorphisms of graph $C^{\ast }$ -algebras.
\newblock {\em Ergodic Theory and Dynamical Systems}, published online 2017 (doi:10.1017/etds.2016.141), 21 pages.

\bibitem[Boyle, 1983]{Boyle1983}
Boyle, M. (1983).
\newblock {\em Topological orbit equivalence and factor maps in symbolic
  dynamics}.
\newblock PhD thesis, University of Washington.

\bibitem[Boyle et~al., 2017]{BoCaEi2017}
Boyle, M., Carlsen, T.~M., \& Eilers, S. (2017).
\newblock Flow equivalence and isotopy for subshifts.
\newblock {\em Dyn. Syst.}, 32(3), 305--325.

\bibitem[Boyle \& Tomiyama, 1998]{BoTo1998a}
Boyle, M. \& Tomiyama, J. (1998).
\newblock Bounded topological orbit equivalence and {$C^*$}-algebras.
\newblock {\em J. Math. Soc. Japan}, 50(2), 317--329.

\bibitem[Brownlowe et~al., 2017]{BrCaWh2017}
Brownlowe, N., Carlsen, T.~M., \& Whittaker, M.~F. (2017).
\newblock Graph algebras and orbit equivalence.
\newblock {\em Ergodic Theory Dynam. Systems}, 37(2), 389--417.

\bibitem[Brix \& Carlsen]{CarlsenBrix}
Brix, K.~A. \& Carlsen, T.~M (2017). 
\newblock Cuntz--Krieger algebras and one-sided conjugacy of shifts of finite type and their groupoids.
\newblock {\em ArXiv e-prints} (arXiv:1712.00179v1), 8 pages.

\bibitem[Carlsen et~al., 2016]{CaEiOr2016}
Carlsen, T.~M., Eilers, S., Ortega, E., \& Restorff, G. (2016).
\newblock {Flow equivalence and orbit equivalence for shifts of finite type and
  isomorphism of their groupoids}.
\newblock {\em ArXiv e-prints} (arXiv:1610.09945v5), 25 pages.

\bibitem[Carlsen \& Rout, 2017]{CaRo2017}
Carlsen, T.~M. \& Rout, J. (2017).
\newblock Diagonal-preserving gauge-invariant isomorphisms of graph
  {$C^\ast$}-algebras.
\newblock {\em J. Funct. Anal.}, 273(9), 2981--2993.

\bibitem[Carlsen et~al., 2017]{CaRuSi2017}
Carlsen, T.~M., Ruiz, E., \& Sims, A. (2017).
\newblock Equivalence and stable isomorphism of groupoids, and
  diagonal-preserving stable isomorphisms of graph {$C^*$}-algebras and
  {L}eavitt path algebras.
\newblock {\em Proc. Amer. Math. Soc.}, 145(4), 1581--1592.

\bibitem[Carlsen et~al.]{CaRuSiTo2017}
Carlsen, T.~M., Ruiz, E., Sims, A., \& Tomforde, M.
\newblock Reconstuction of groupoids and $C^*$-rigidity of dynamical systems.
\newblock {\em ArXiv e-prints} (arXiv:1711.01052v1), 44 pages.

\bibitem[Carlsen \& Winger, 2016]{CaWi2016}
Carlsen, T.~M. \& Winger, M.~L. (2016).
\newblock {Orbit equivalence of graphs and isomorphism of graph groupoids}.
\newblock {\em ArXiv e-prints} (arXiv:1610.09942v3), 9 pages.

\bibitem[Cuntz, 1981]{Cu1981}
Cuntz, J. (1981).
\newblock A class of {$C\sp{\ast} $}-algebras and topological {M}arkov chains.
  {II}. {R}educible chains and the {E}xt-functor for {$C\sp{\ast} $}-algebras.
\newblock {\em Invent. Math.}, 63(1), 25--40.

\bibitem[Cuntz \& Krieger, 1980]{CuKr1980}
Cuntz, J. \& Krieger, W. (1980).
\newblock A class of {$C\sp{\ast} $}-algebras and topological {M}arkov chains.
\newblock {\em Invent. Math.}, 56(3), 251--268.

\bibitem[Giordano et~al., 1995]{GiPuSk1995}
Giordano, T., Putnam, I.~F., \& Skau, C.~F. (1995).
\newblock Topological orbit equivalence and {$C\sp *$}-crossed products.
\newblock {\em J. Reine Angew. Math.}, 469, 51--111.

\bibitem[Hahn, 1978]{Ha1978}
Hahn, P. (1978).
\newblock Haar measure for measure groupoids.
\newblock {\em Trans. Amer. Math. Soc.}, 242, 1--33.

\bibitem[an~Huef \& Raeburn, 1997]{HuRa1997}
an~Huef, A. \& Raeburn, I. (1997).
\newblock The ideal structure of {C}untz-{K}rieger algebras.
\newblock {\em Ergodic Theory Dynam. Systems}, 17(3), 611--624.

\bibitem[Kitchens, 1998]{Ki1998}
Kitchens, B.~P. (1998).
\newblock {\em Symbolic dynamics}.
\newblock Universitext. Berlin: Springer-Verlag.
\newblock One-sided, two-sided and countable state Markov shifts.

\bibitem[Kumjian, 1986]{Ku1986a}
Kumjian, A. (1986).
\newblock On {$C^\ast$}-diagonals.
\newblock {\em Canad. J. Math.}, 38(4), 969--1008.

\bibitem[Li, 2016]{Li2016}
Li, X. (2016).
\newblock Continuous orbit equivalence rigidity.
\newblock {\em Ergodic Theory and Dynamical Systems}, published online 2016 (doi:10.1017/etds.2016.98), 21 pages.

\bibitem[Lind \& Marcus, 1995]{LiMa1995}
Lind, D. \& Marcus, B. (1995).
\newblock {\em An introduction to symbolic dynamics and coding}.
\newblock Cambridge: Cambridge University Press.

\bibitem[Matsumoto, 2010]{Ma2010}
Matsumoto, K. (2010).
\newblock Orbit equivalence of topological {M}arkov shifts and
  {C}untz-{K}rieger algebras.
\newblock {\em Pacific J. Math.}, 246(1), 199--225.

\bibitem[Matsumoto, 2016]{Ma2016c}
Matsumoto, K. (2016).
\newblock Continuous orbit equivalence, flow equivalence of Markov shifts and
  circle actions on Cuntz--Krieger algebras.
\newblock {\em Mathematische Zeitschrift}, 285(1--2), 121--141.

\bibitem[Matsumoto \& Matui, 2016]{MaMa2016}
Matsumoto, K. \& Matui, H. (2016).
\newblock Continuous orbit equivalence of topological {M}arkov shifts and
  dynamical zeta functions.
\newblock {\em Ergodic Theory Dynam. Systems}, 36(5), 1557--1581.

\bibitem[Parry \& Sullivan, 1975]{PaSu1975}
Parry, B. \& Sullivan, D. (1975).
\newblock A topological invariant of flows on {$1$}-dimensional spaces.
\newblock {\em Topology}, 14(4), 297--299.

\bibitem[Paterson, 1999]{Pa1999}
Paterson, A. L.~T. (1999).
\newblock {\em Groupoids, inverse semigroups, and their operator algebras},
  volume 170 of {\em Progress in Mathematics}.
\newblock Boston, MA: Birkh\"auser Boston Inc.

\bibitem[Quigg, 1996]{Qu1996}
Quigg, J.~C. (1996).
\newblock Discrete {$C\sp *$}-coactions and {$C\sp *$}-algebraic bundles.
\newblock {\em J. Austral. Math. Soc. Ser. A}, 60(2), 204--221.

\bibitem[Renault, 1980]{Re1980a}
Renault, J. (1980).
\newblock {\em A groupoid approach to {$C\sp{\ast} $}-algebras}, volume 793 of
  {\em Lecture Notes in Mathematics}.
\newblock Berlin: Springer.

\bibitem[Renault, 2008]{Re2008}
Renault, J. (2008).
\newblock Cartan subalgebras in {$C\sp *$}-algebras.
\newblock {\em Irish Math. Soc. Bull.}, (61), 29--63.

\bibitem[Sims, 2017]{Sims2017}
Sims, A. (2017).
\newblock Hausdorff {\'e}tale groupoids and their $C^*$-algebras.
\newblock {\em To appear in the volume ``Operator algebras and dynamics:
  groupoids, crossed products and Rokhlin dimension'' in Advanced Courses in
  Mathematics. CRM Barcelona.}

\bibitem[Tomiyama, 1996]{To1996}
Tomiyama, J. (1996).
\newblock {$C^*$}-algebras and topological dynamical systems.
\newblock {\em Rev. Math. Phys.}, 8(5), 741--760.

\end{thebibliography}

\end{document}